\documentclass[a4paper]{article}
\usepackage{amsmath} % because of {split}. amssymb is not needed now?
\usepackage{amsfonts}
\def\X{\Xi^1}
\def\T{\Xi^2}
\def\U{\Phi}
\def\Ui{\Phi^1}
\def\Real{\mathbb{R}}
\def\grad#1#2{\mathop{\rm grad}\nolimits_{#1}#2}
\def\grad#1#2{\frac{\partial #2}{\partial #1}}
\newenvironment{longmath}{
  \begin{center} \begin{math} \displaystyle
}{
  \end{math} \end{center}
}
\def\supp{\mathop{\rm supp}\nolimits}
\def\id{\mathop{\rm id}\nolimits}

\def\prolong^#1{\mathop{\rm pr}\nolimits^{(#1)}}
\def\dual<#1,#2>{\left\langle #1,\, #2 \right\rangle}
\def\class[#1]{\mathop{\rm cl}\left[ #1 \right]}
\def\abs|#1|{\left| #1 \right|}
\def\Dj{\rlap{--}D}

\def\epsilon{\varepsilon}
\usepackage{latexsym} 
\usepackage{amssymb} 
\newcommand{\ep}{\hspace*{\fill}$\Box$} 
\newcommand{\eps}{\varepsilon} 
\newcommand{\pr}{{\bf Proof. }} 
\newcommand{\ms}{\medskip\\} 
\newcommand{\cl}{\mbox{\rm cl}}

\newcommand{\prol}{\mbox{\rm pr}^{(n)}}

\newcommand{\R}{\mathbb{R}}
\newcommand{\N}{\mathbb{N}} 
\newcommand{\C}{\mathbb{C}}

\newcommand{\be}{\begin{equation}} 
\newcommand{\ee}{\end{equation}} 
\newcommand{\bea}{\begin{eqnarray}} 
\newcommand{\eea}{\end{eqnarray}} 
\newcommand{\beast}{\begin{eqnarray*}} 
\newcommand{\eeast}{\end{eqnarray*}}

\newtheorem{thr}{\hspace*{-1.1mm}}[section] 
\newcommand{\bt}{\begin{thr} {\bf Theorem }} 
\newcommand{\et}{\end{thr}} 
\newcommand{\bp}{\begin{thr} {\bf Proposition }} 
\newcommand{\bc}{\begin{thr} {\bf Corollary }} 
\newcommand{\blem}{\begin{thr} {\bf Lemma }} 
\newcommand{\bex}{\begin{thr} {\bf Example }\rm} 
\newcommand{\bexs}{\begin{thr} {\bf Examples }\rm} 
\newcommand{\brm}{\begin{thr} {\bf Remark }\rm} 
\newcommand{\brms}{\begin{thr} {\bf Remarks }\rm} 
\newcommand{\bd}{\begin{thr} {\bf Definition }} 
\newcommand{\Om}{\Omega} 
\newcommand{\D}{{\cal D}}
\newcommand{\G}{{\cal G}}
\newcommand{\comp}{\subset\subset}
\newcommand{\vphi}{\varphi}

\newcommand{\eigmat}{{\mathrm I}}% Identity matrix
\newcommand{\lie}[1]{\frac{\partial}{\partial \eta} \Big|_{\eta=0} #1}

\newcommand{\dimension}[3]{[#1]^{#2}_{#3}}

\title{Symmetry group analysis of weak solutions}
\author
{N.\ \Dj api\'c \and
 M.\ Kunzinger  \and  
 S.\ Pilipovi\'c}
\date{%\today; preliminary version
}

\def\to{\rightarrow}

\newcounter{bean} % bean is from Lamport's example of list

\begin{document}
\maketitle

\begingroup
\renewcommand{\thefootnote}{} % do not print footnote numbers
% N.\ \Dj api\'c
% \refstepcounter{footnote}
 \footnotetext{Institute~of~Mathematics,
  University~of~Novi Sad,
  Trg~D.~Obradovi\'ca~4,
  21000~Novi~Sad,
  Yugoslavia.
  E-mail: {\tt nenad@unsim.ns.ac.yu}} \and
% M.\ Kunzinger
% \refstepcounter{footnote}
 \footnotetext{Institut f\"ur Mathematik,
  Universit\"at Wien, 
  Strudlhofg.~4,
  A-1090 Wien,
  Austria.
  E-mail: {\tt Michael.Kunzinger@univie.ac.at}} \and  
% S.\ Pilipovi\'c
% \refstepcounter{footnote}
 \footnotetext{Institute~of~Mathematics,
  University~of~Novi Sad,
  Trg~D.~Obradovi\'ca~4,
  21000~Novi~Sad,
  Yugoslavia.
  E-mail: {\tt pilipovic@unsim.ns.ac.yu}}
\endgroup % now print footnote numbers

\begin{abstract}
Methods of Lie group analysis of differential equations are extended to weak
solutions of (linear and nonlinear) PDEs, where the term ``weak solution'' comprises
the following settings:
\begin{itemize}
\item[(a)] Distributional solutions.
\item[(b)] Solutions in  generalized function algebras.
\item[(c)] Solutions in the sense of
association (corresponding to a number of weak or integral
solution concepts in classical analysis).
\end{itemize}
Factorization properties and infinitesimal criteria are developed that allow to treat all
three settings simultaneously, thereby unifying and extending previous work in this area.
\end{abstract}
\section{Introduction} \label{intro}

Local symmetries for equations with weak type solutions, such as, e.g., 
conservation laws,
involve different constraints depending on the framework
in which the equations are analyzed. The aim of this paper is to study symmetry 
properties of differential equations involving singular 
(in particular: distributional) 
terms through an analysis of symmetries in distribution spaces and generalized 
function algebras, as well as associated 
(i.e., weak type) symmetries. 
A main ingredient in our analysis will be
the determination of infinitesimal criteria for these solution
concepts.

Investigations in this
direction have been initiated by Meth\'{e} (\cite{Me}), Tengstrand (\cite{T}), 
Szmydt and Ziemian (\cite{S,S2,SZ}) and have been systematically pursued by 
Berest and Ibragimov (\cite{huy,wi,ga,bi,amp, Ibragimov}) in the distributional 
setting. More recently, in \cite{Kunzinger, Kunzinger-O, KO}, 
an extension of the purely distributional 
methods applied so far has been given that allows to also consider nonlinear 
equations involving singularities.
The basic tool allowing for such an extension
is Colombeau's theory of algebras of generalized functions.
In what follows, on the one hand we are going to continue the analysis of 
\cite{KO}, and on the other hand we shall establish connections between
the distributional criteria developed in \cite{wi,ga} and the Colombeau-type 
methods given in \cite{KO}.  Moreover, the
%Colombeau
framework of generalized function algebras enables us to study
symmetry properties of associated solutions (i.e.\ of weak or integral solutions) 
by the same methods.

To begin with, let us fix some notations concerning group analysis of differential
equations and Colombeau's theory of generalized functions. Our principal reference
for symmetries of differential equations is \cite{Olver}, whose terminology we
shall follow closely. Let
\(\cal M \) be an open subset of \( \Real ^{p+q} \) (in what follows, $p$ will be
the number of independent variables of a system of differential equations, $q$ the
number of dependent variables) and  \( G \) a Lie group acting regularly on $\cal M$.

For \( x \in \Real^p \), \( u \in \Real^q \) and \( (x,u) \in {\cal M} \) we write
\[
  (\tilde x, \tilde u) =
  g \cdot (x,u) = (\Xi_g (x,u), \Phi_g (x,u))
\]
If \(\Xi_g (x,u) =   \Xi_g (x) \), \( g \in G \), then \( G \) is called projectable.
Elements of the Lie algebra of $G$ as well as the corresponding vector 
fields   on   $\cal M$   will  typically be  denoted  by  ${\bf  v}$.

Put \( {\cal X} = \Real^p \), \( {\cal U} = \Real^q \). By identifying a function
\( f : \Omega \subset {\cal X} \to {\cal U} \) with its graph
\(  \Gamma _f = \{(x, f (x)) \mid x \in \Omega \} \subset {\cal X} \times {\cal U}\),
the action of $g \in G$ onto $f$ is defined (locally) by
\[
  g \cdot f = \tilde f = \left(
    \Phi_g \circ (\id _{\cal X} \times f)
  \right) \circ \left(
    \Xi_g \circ (\id _{\cal X} \times f)
  \right) ^{-1}
\]
where
\( \id _{\cal X} \) is the identity mapping on
\( {\cal X} \); in the projectable case this specializes to
\[
  g \cdot f =
  \Phi_g \circ
  (\id _{\cal X} \times f) \circ
  \Xi_g ^{-1}.
\]

Set \({\cal U} ^{(n)} = {\cal U} \times {\cal U} _1 \times \ldots \times {\cal U} _n
\), \({\cal M} ^{(n)} =   {\cal M} \times {\cal U} _1 \times \ldots \times 
{\cal U} _n \).
%, where \( {\cal M} \) is open in \( {\cal X} \times {\cal U} \).
Here \( {\cal U} _k = \Real ^{q p_k} \), with coordinates \( u_J^\alpha \),
\( \alpha = 1,\ldots,q \), \(J = (j_1,\ldots,j_k) \in \{1,\ldots,p\} ^k
\) such that \( j_1 \leq j_2 \leq \ldots \leq j_k \) and
\( p_k = {p+k-1 \choose k} \), is the number of different partial
derivatives of order \( k \) of a scalar valued smooth function of
\( p \) variables. Elements of \( {\cal U} ^{(n)} \) are denoted by \( u^{(n)} \). 
${\cal M} ^{(n)}$ is called the $n$-jet space of $\cal M$ 
and we set $N:=\mbox{dim}({\cal M}^{(n)}).$
The coordinates on ${\cal M}^{(n)}$ will also be written as 
$(z_1,\dots,z_p,z_{p+1},\dots,z_N) = (x,u^{(n)})$.
For any $f: \Omega \subseteq {\cal X} \to \R^q$, the $n$-th prolongation
or $n$-jet of $f$ is the function $\prolong^n f:\Omega \to {\cal U}^{(n)}$ formed
by $f$ and its derivatives up to order $n$. The $n$-th prolongation of a 
group action $g$ or vector field ${\bf  v}$ is written as  
$\prol g$ or $\prol {\bf  v}$, respectively.

Let \( \Sigma \) be a system of differential equations with
\( p \) variables and \( q \) unknown functions of the form
\begin{equation} \label{system}
 \Delta_i (x,u^{(n)}) = 0, \quad
 i=1,\ldots,s
\end{equation}
(with $\Delta_i$ smooth for all $i$). We shall henceforth assume that
(\ref{system}) is nondegenerate (i.e. locally solvable and of maximal rank,
see \cite{Olver}).
Any \( f : \Omega \subset {\cal X} \to {\cal U} \) which solves the
system on its domain will be called a solution. This amounts to saying that
the graph of the $n$-jet of $f$, $\Gamma_f^{(n)}$ is contained in the zero-set
$\Sigma_\Delta$ of $\Delta$. A symmetry group of \( \Sigma \) is a local 
transformation group on \( {\cal X} \times {\cal U} \) such that if \( f \) 
is a solution of the system, \( g \in G \) and \( g \cdot f \) is defined then 
also \( g \cdot f \) is a solution of \( \Sigma \).

Next, let us shortly recall some basic definitions from Colombeau's theory of
generalized functions (\cite{AB}, \cite{Biagioni}, \cite{c1}, \cite{Colombeau:E}, \cite{c3},
\cite{Oberguggenberger}). For notational simplicity we are going to work in the so-called 
``special'' Colombeau algebra ${\cal G}(\Omega)$, defined as the quotient algebra
${\cal E}_M(\Omega)/{\cal N}(\Omega)$, where 
\begin{eqnarray*}  
{\cal E}_M(\Omega)&:=&\{(u_\varepsilon)_{\eps\in I}\in 
{\cal C}^\infty(\Omega)^I: \forall K\subset\subset\Omega,  
\forall\alpha\in\N_o^n 
\mbox{ }\exists p\in \N \mbox{ with }\\&&\sup_{x\in K}|\partial^\alpha  
u_\varepsilon(x)|=O 
(\varepsilon^{-p})\mbox{ as }\varepsilon\rightarrow 0\}\\{\cal N}(\Omega)&:=& 
\{(u_\varepsilon)_{\eps\in I}\in {\cal C}^\infty(\Omega)^I:  
\forall K\subset\subset\Omega, \forall\alpha\in\N_o^n\mbox{ }\forall q\in\N\\  
&&\sup_{x\in K}|\partial^\alpha u_\varepsilon(x)|=O(\varepsilon^{q})\mbox{ as } 
\varepsilon\rightarrow 0\}
\end{eqnarray*} 
Here $I = (0,1]$. ${\cal G}(\Omega)$ is a differential algebra (with componentwise
operations) and $\Omega \to {\cal G}(\Omega)$ is a fine sheaf on $\R^n$. The
equivalence class of $(u_\eps)_{\eps\in I}$ in ${\cal G}(\Omega)$ will be denoted
by $U=\cl[(u_\eps)_{\eps\in I}]$ or $\cl[u_\eps]$ for short. 

We shall make use of the function spaces 
\begin{eqnarray*}
%\begin{array}{l}
{\cal D}(\Omega) &=& \{f\in {\cal C}^\infty(\Omega): 
\supp(f) \ \mbox{compact}\} 
\\
{\cal S}(\R^n) &=& \{f\in {\cal C}^\infty(\R^n): 
\forall q>0 \ \forall \alpha\in
\N_0^n\ \sup\limits_{x\in \R^n} (1+|x|)^q 
|\partial^\alpha f(x)|<\infty\} \\
{\cal O}_M(\Omega) &=& \{f\in {\cal C}^\infty(\Omega):
\forall \alpha\in \N_0^n \ \exists p>0
\ \sup\limits_{x\in \Omega} (1+|x|)^{-p} 
|\partial^\alpha f(x)|<\infty\} 
%\end{array}  
\end{eqnarray*}
$\D(\Omega)$ is the space of test functions on $\Om$, elements of ${\cal S}(\R^n)$
and ${\cal O}_M(\R^n)$ are called {\em rapidly decreasing} and
{\em slowly increasing}, respectively.

%Furthermore, we shall need the algebra  
The algebra 
${\cal G}_\tau(\Omega)={\cal E}_\tau(\Omega)/{\cal N}_\tau(\Omega)$ 
of tempered generalized functions is defined by
\begin{eqnarray*} 
&& {\cal E}_\tau(\Omega) = \{ (u_\eps)_{\eps\in I}\in  
({\cal O}_M(\Omega))^I : 
\forall \alpha\in \N_o^n\ \exists p>0 \\ 
&&\qquad\qquad\quad\sup_{\scriptstyle x\in \Omega}(1+|x|) 
^{-p}|\partial^\alpha u_\eps(x)| = O(\eps^{-p})\ (\eps\to 0)\}\\ 
&& {\cal N}_\tau(\Omega) = \{ (u_\eps)_{\eps\in I}\in  
({\cal O}_M(\Omega))^I : 
\forall \alpha\in \N_o^n\ \exists p>0\ \forall \ q>0 \\ 
&&\qquad\qquad\quad\sup_{\scriptstyle x\in \Omega}(1+|x|) 
^{-p}|\partial^\alpha u_\eps(x)| = O(\eps^{q})\ (\eps\to 0)\}
\end{eqnarray*}

By componentwise insertion, elements of $\G(\Om)$ and $\G_\tau(\Om)$
can be composed  with slowly increasing functions. 
Next, choose $\rho \in {\cal S}(\R^n)$ 
such that $\int \rho(x) \, dx = 1$ and
$\int \rho(x) x^\alpha\, dx = 0$ for all $\alpha \in \N_0^n$ with $|\alpha|\ge 1$.
Then ${\cal E}'(\Omega)$ (the space of compactly supported distributions)
is linearly embedded into ${\cal G}(\Om)$ via   
$\iota:   u   \to  \cl[(u\ast\rho_\eps)_{\eps\in I}]$ (where $\rho_\eps = \frac{1}{\eps^n} 
\rho(\frac{.}{\eps})$). Moreover $\iota$ coincides with the identical embedding $\sigma:
f \to \cl[(f)_{\eps\in I}]$ on ${\cal D}(\Omega)$, so ${\cal D}(\Omega)$
becomes a subalgebra of ${\cal G}(\Omega)$ via $\iota$. Finally, there is a unique
sheaf morphism $\hat \iota$ extending $\iota$
to ${\cal C}^\infty(\,.\,) \hookrightarrow {\cal D}'(\,.\,) \hookrightarrow 
{\cal G}(\,.\,)$ (where $\D'$ denotes the space of distributions). $\hat \iota$ commutes 
with partial derivatives, and its restriction to ${\cal C}^\infty$ is a sheaf morphism of algebras.
The map $\iota$ defined above also provides a linear embedding  
of ${\cal S}'(\R^n)$ into ${\cal G}_\tau(\R^n)$  
commuting with partial derivatives and making 
\[ 
{\cal O}_C(\R^n) = \{f\in {\cal C}^\infty(\R^n): 
\exists p>0 \ \forall \alpha\in 
\N_o^n\ \sup\limits_{\scriptstyle x\in \R^n}  
(1+|x|)^{-p} |\partial^\alpha f(x)|<\infty\} 
\] 
a faithful subalgebra.
So far, Colombeau algebras are the only known differential algebras enjoying these
optimal embedding properties. Moreover, an intrinsic global formulation of Colombeau's
construction on differentiable manifolds retaining all the characteristics of the
local theory has recently been achieved (\cite{found,vim}). 
Of the various variants of Colombeau algebras we shall
also use the subalgebra ${\cal G} _{\infty} (\Omega)$ of $\G(\Om)$ consisting of
those elements of $\G(\Om)$ possessing a representative $(u_\eps)_\eps$ such that
$\sup_{\eps\in I}\|u_\eps\|_{L^\infty(\Om)} < \infty$. Finally, we shall 
employ the ``mixed type''-algebra
$\tilde \G(\R\times\Om)$ whose elements satisfy $\G$-bounds with respect
to $t\in \R$ and $\G_\tau$-bounds with respect to $x\in \Om$.

The ring of constants in ${\cal G}$ is denoted by $\overline{\C}$, its elements are
called generalized numbers. Elements of $\overline{\C}$ may be used to model
infinitesimal numbers (e.g., $(\eps)_{\eps\in I}$ is a representative of an infinitely
small yet nonzero generalized number), which may be viewed as a 
``nonstandard'' aspect of the theory.
The support of  a generalized function $U\in {\cal G}(\Om)$, 
$\supp_g U$, is defined as the complement of the largest open
subset $\Omega'$ of $\Om$  such that $U_{|\Omega'}=0$.
This notion is coherent with the embedding $\iota$, i.e. for any $T\in{\cal D}'(\Omega)$ we 
have $\supp T = \supp_g (\iota(T))$. Finally, we mention the concept of association in
the algebra $\cal G$: Two elements $F$, $G$ are associated
$F \approx G$ if there exist representatives $F_\varepsilon$ and
$G_\varepsilon$ of $F$ and $G$, respectively, such that
$$
  \lim _{\varepsilon \to 0}
  \int (G_\varepsilon (x) - F_\varepsilon (x)) \psi (x)dx = 0,
  \mbox{ for all } \psi\in {\cal D} (\Omega)
$$
Clearly, this definition does not depend on the choice of
representatives. The concept of association (resp.\ strong association,
cf.\ \ref{strass} below) plays a central role in the Colombeau framework as in many
cases it allows for a distributional interpretation of results achieved in $\G$.
Particularly in applications to physics and numerics it is of fundamental
importance (cf.\ e.g.\ \cite{Biagioni}, \cite{c3}, \cite{V}).

With this terminology at hand we can now formulate the main goals of this article:
Let $G$ be a transformation group acting on the space of independent and dependent 
variables of a system (\ref{system}) of (linear or nonlinear)
differential equations. We are looking for criteria for $G$ to transform weak
solutions of (\ref{system}) into weak solutions. More precisely, we shall develop
conditions under which $G$ transforms $\D'$-solutions to $\D'$-solutions (in
case (\ref{system}) is linear), $\G$-solutions to $\G$-solutions, or solutions
in the sense of association into solutions in the sense of association (in which 
case (\ref{system}) is to be replaced by $\Delta(x,u^{(n)}) \approx 0$). We shall
see that these questions are in fact closely linked and that criteria for one
situation are often useful in other cases as well.
As an application we study the symmetries of the quasilinear hyperbolic
system $u_t+A(u)u_x=0$ where $A$ is an $s\times s$ matrix with ${\cal C}^1(\Real^s)$
entries (cf.\ \cite{Smol}, \cite{Bres}). 
For treating this system distribution type spaces are not convenient,
while ${\cal G} $ provides a quite  satisfactory solution concept
(cf.\ \cite{Ke}, \cite{Ned} and the literature cited therein). We consider
strongly associated solutions and calculate symmetry transformations of such
solutions. Also, we discuss infinitesimal criteria for symmetries of this system.

\section{Factorization properties, symmetries in ${\cal G}$}

Let $G$ be a projectable group action on some open set ${\cal M} \subset {\cal X}\times 
{\cal U}$. As was already pointed out in the introduction,
composition in the framework of Colombeau generalized functions
requires polynomial growth restrictions on the smooth functions. 
Thus we first single out those group actions which can be applied to elements
of $\G$ (cf.\ \cite{KO}). An element $g \in G$ with
\(
  g \cdot (x,u) = (\Xi_g (x), \Phi_g (x,u))
\),
is called slowly increasing if $u \mapsto \Phi_g (x,u)$ is slowly increasing,
uniformly for $x$ in compact sets; $g$ is called strictly slowly increasing if
$\Phi _g \in {\cal O} _M ({\cal M})$.

If $\Omega\subset {\cal X}$, $U\in {\cal G}(\Omega)$ and $g$ is slowly increasing,
the action of $g$ on $U$ is the element of ${\cal G}(\Xi_g(\Omega))$ given by
\begin{equation} \label{gafuncprgen}
  g U = \class[
    \left(
      \Phi_g \circ
      (\id \times u_\varepsilon) \circ
      \Xi_g^{-1}
    \right) _\varepsilon
  ]
\end{equation}

Also, in order to be able to insert elements of $\G(\Om)$ into (\ref{system}) we
will from now on suppose that $u^{(n)} \mapsto \Delta(x,u^{(n)})$ is slowly increasing,
uniformly for $x$ in compact sets. A symmetry group of (\ref{system}) in ${\cal G}$
is a local transformation group acting on ${\cal X}\times {\cal U}$ such that
if $U$ is a solution of the system in ${\cal G}$, $g\in G$ and $g\cdot U$ is
defined, then also $g\cdot U$ is a solution of $\Sigma$ in ${\cal G}$.

Let $G$ be a slowly increasing symmetry group of some
differential equation
\begin{equation} \label{diffe}
  \Delta(x,u^{(n)}) = 0
\end{equation}
and let $U\in {\cal G}(\Omega)$ be a generalized solution to (\ref{diffe}).
Then for any representative
$(u_\varepsilon)_\varepsilon$ of $U$ there exists some
$(n_\varepsilon)_\varepsilon \in {\cal N}(\Omega)$ such that for all
$x$ and all $\varepsilon$ we have
\begin{equation} \label{differep}
  \Delta(x,\prolong^n u_\varepsilon(x)) = n_\varepsilon(x)
\end{equation}
Due to the nontrivial right hand side of (\ref{differep}) it is clear that a 
direct (componentwise) application of classical symmetry methods to Colombeau 
solutions is not feasible. A transfer of classical symmetry groups into the
$\G$-setting therefore has to rely on properties of symmetry transformations
that are better suited to the algebraic structure of $\G$. The key concept 
serving this purpose (and, as we shall see shortly, at the same time 
applicable to symmetries of $\D'$- and other weak solutions)
is that of factorization:

Let $G = \{ g_{\eta} \mid \eta \in \Real \}$
be a (classical) one parameter symmetry group of (\ref{system})
\[
 g_\eta\cdot (x,u) = (\Xi_\eta(x,u),\Phi_\eta(x,u))
\] 
Then by \cite{Olver}, Prop.\ 2.10  there exists a smooth map $Q:
{\cal W} \to \R^{s^2}$ (${\cal W}$ open in $\R\times {\cal M}^{(n)}$,
$0\times {\cal M}^{(n)} \subseteq {\cal W}$) such that 
\begin{equation}\label{factor0}
\Delta(\prol g_\eta(z)) = Q(\eta,z)\Delta(z) 
\end{equation}
Throughout this paper, ${\cal W}$ will denote an open set as specified above.
By \cite{KO}, Th.\ 3.4 for any smooth $u : \Omega \subset \Real ^p\to \Real ^q$ 
such that $g_\eta u$ exists we have
\begin{equation}\label{funcfactor}
\begin{array}{l}
\Delta(\Xi_\eta(x,u(x)),\prol(g_\eta u)(\Xi_\eta(x,u(x)))) \\
=Q(\eta,x,\prol u(x)) \cdot \Delta(x,\prol u(x))\,.  
\end{array}
\end{equation}
%for $1 \le \nu \le s$. 
In particular, $\eta\to g_\eta$
is a symmetry group of (\ref{system}) in $\G$ if $(x,u^{(n)}) \to 
Q_{\mu \nu}(\eta,x,$ $\prolong^n u(x))$ is slowly increasing with respect to $u^{(n)}$,
uniformly for $x$ varying in compact sets (\cite{KO}, Proposition 3.5). Thus the
need for conditions ensuring that the $Q_{\mu\nu}$ remain well behaved 
(in the above sense)
arises. Theorem 3.8 of \cite{KO} shows that for scalar differential equations possessing
a stand alone term (i.e. $\frac{\partial\Delta}{\partial z_k} = c$ for some $k>p$)
this is indeed always the case. Our first aim is to generalize this result.

Since (\ref{system}) is nondegenerate it follows that the Jacobian
$J(\Delta)$ of $\Delta$ has rank $s$ on the zero-set of $\Delta$. The following 
result uses a mild strengthening of this assumption to derive a factorization
property adapted to the polynomial growth restrictions necessary for
applying nonlinearities to elements of $\G$.

%%%%%%%%%%%%%%%%%%%%%%%%%%%%%%%%%%%%%%%%%%%%%%%%%%%%%%%%%%%%%%%%%%%%%%%%%%
\bt \label{factorth}
Let   $G = \{ g_{\eta } \mid \eta \in (-\eta _{0} , \eta _{0})\}$ be a 
slowly increasing symmetry group of system (\ref{system}) and suppose that
there exist $p < k_1< \dots < k_s \le N$ such that, setting
$J_{k_1,\dots,k_s}(\Delta) := \left(\frac{\partial 
\Delta_i}{\partial z_{k_j}}\right)_{i,j=1\dots s}$, the following conditions are 
satisfied:
\begin{itemize}
\item[(i)] $z \to (z',\Delta(z))$ is injective, where $z' = (z_1,\dots,\widehat{z_{k_1}},
\dots,\widehat{z_{k_s}},\dots,z_N)$.
\item[(ii)] $z \to (\det(J_{k_1,\dots,k_s}(\Delta)))^{-1} (z)$
is defined globally and is slowly 
increasing, uniformly for $(z_1,\dots,z_p)$ varying in compact sets.
\end{itemize}
Then there exists a smooth mapping $Q: {\cal W} \to \R^{s^2}$ 
which is slowly increasing in $z\in {\cal M}^{(n)}$, uniformly for $z_1,\dots,z_p$
varying in compact sets such that (\ref{factor0}) holds. In particular, 
(\ref{funcfactor}) holds for any smooth $u: \Om\subseteq \R^p \to \R^q$ such that
$g_\eta u$ exists.
\et
%%%%%%%%%%%%%%%%%%%%%%%%%%%%%%%%%%%%%%%%%%%%%%%%%%%%%%%%%%%%%%%%%%%%%%%%%%
\pr Without loss of generality we may suppose $\{k_1,\dots,k_s\} = 
\{N-s+1,\dots,N\}$. We set $z' = (z_1,\dots,z_{N-s})$,
$z'' = (z_{N-s+1},$ $\dots,$ $z_N)$ and define $\tilde \Delta: \R^N \to \R^N$ by
$$
\tilde \Delta(z) = (z',\Delta(z))
$$
Since $\det(J(\tilde\Delta)) = \det(J_{k_1,\dots,k_s}(\Delta)) \not=0$, $\tilde \Delta$
is a diffeomorphism by (i). Moreover, 
$$
\Delta\circ \tilde\Delta^{-1} = (y_1,\dots,y_N) \to (y_{N-s+1},\dots,y_N) = y \to y''
$$
Now set $f_\eta(z) = \Delta(\prol g_\eta(z))$.
By \cite{KO}, Proposition 3.7 and our general assumption on $\Delta$, $f_\eta$ 
is slowly increasing in $z$, uniformly for $(z_1,\dots,z_p)$ varying in compact sets.
Since $\eta \to g_\eta$ is a symmetry group of (\ref{system}) we have
$f_\eta \circ \tilde\Delta^{-1}(y) \equiv 0$ if $y''\equiv 0$. Thus
\begin{equation}\label{f}
\begin{array}{rcl}
f_\eta\circ\tilde\Delta^{-1}(y) &=& (f_\eta\circ\tilde\Delta^{-1})(y',\tau y'')
\mid_{\tau=0}^1 \\
&=& \int_0^1 \frac{d}{d\tau} (f_\eta\circ\tilde\Delta^{-1})(y',\tau y'')\,d\tau \\
&=& \left(\int_0^1 J_{k_1\dots k_s}(f_\eta\circ\tilde\Delta^{-1})(y',\tau y'')
\,d\tau\right)\cdot y''
\end{array}
\end{equation}
Inserting $z = \tilde\Delta^{-1}(y)$ into (\ref{f}) we arrive at
\begin{equation}\label{f2}
f_\eta(z) = \left(\int_0^1 J_{k_1,\dots, k_s}(f_\eta\circ\tilde\Delta^{-1})
(z',\tau \Delta(z))\,d\tau\right)\cdot \Delta(z) =: Q(\eta,z)\cdot \Delta(z)
\end{equation}
In order to establish the claimed growth properties of $Q$, since $f_\eta$ and
$\Delta$ are slowly increasing, by the chain rule it suffices to estimate
$J_{k_1\dots k_s}(\tilde\Delta^{-1})$. The determinant of this Jacobian is
precisely $(\det(J_{k_1,\dots,k_s}(\Delta)))^{-1}$, so the claim follows from
(ii). \ep 

\brms \label{rem1}
\begin{itemize}
\item[(i)] By dropping the growth restrictions on $\Delta$, $G$ (in particular,
allowing for $G$ to be nonprojectable and merely supposing that
$(\det(J_{k_1,\dots,k_s}(\Delta)))^{-1}$ exists globally), the same proof as above still
provides the explicit form (\ref{f2}) of the factorization property of 
general symmetry groups of (\ref{system}), which reads
\begin{equation}\label{fgeneral}
\begin{array}{l}
\!\!\!\!\!\Delta(\prol g_\eta(z)) = 
(\int_0^1 J_{k_1,\dots, k_s}(f_\eta\circ\tilde\Delta^{-1})
(z_1,\dots, z_{k_1-1},\tau \Delta_1(z),z_{k_1+1},  \\
\!\!\!\!\! 
\qquad \qquad \qquad \qquad
\dots,z_{k_s-1},\tau \Delta_s(z),z_{k_s+1},\dots,z_N)\,d\tau) \cdot \Delta(z) 
\end{array}
\end{equation}
\item[(ii)] A necessary and sufficient condition for (ii) in \ref{factorth} is given by
\begin{equation}\label{a3} 
\begin{array}{l}
\forall K \comp \R^p \ \exists C > 0 \ \exists r >0 \ \forall z\in {\cal M}^{(n)}: \\
\inf_{(z_{1},\dots,z_{p})\in K} \left| 
\frac{\partial(\Delta_1,\dots,\Delta_s)}{\partial(z_{k_1},
\dots,z_{k_s})}(z)\right| \geq C ((1+ |z_{p+1}|) \dots (1+|z_{N}|))^{-r} 
\end{array}
\end{equation}

\item[(iii)] An extensive compilation of sufficient conditions for 
global injectivity of smooth maps ($\tilde \Delta$ in our case) can be found 
in \cite{Part}. As an example we mention a result of Gale and Nikaido
(\cite{Part}, ch.\ 3) stating that any $F: \Om \to \R^n$ 
(with $\Om$ a rectangular region
in $\R^n$) is injective if its Jacobian is a $P$-matrix on all of $\Om$
(i.e.\ all principal minors are positive). Moreover, global invertibility results
for Sobolev functions are given in \cite{Ball}.
\end{itemize}
\et

\bc \label{symmG}
Under the assumptions of Theorem \ref{factorth}, $G$ is a symmetry group 
of (\ref{system}) in $\G$.
\et
%%%%%%%%%%%%%%%%%%%%%%%%%%%%%%%%%%%%%%%%%%%%%%%%%%%%%%%%%%%%%%%%%%%%
\pr
Let $u \in \G(\Om)^q$ be a solution of (\ref{system}). Then by (\ref{funcfactor}),
for any representative $(u_\eps)_{\eps\in I}$ of $u$ we have, taking into account
the projectability of $g_\eta$,
\begin{equation} \label{factorpro}
\begin{array}{l}
\Delta(x,\prol(g_\eta u_\eps)(x))\\ 
=Q(\eta,\Xi_{-\eta}(x),\prol u_\eps(\Xi_{-\eta}(x))) \cdot 
\Delta(\Xi_{-\eta}(x),\prol u_\eps(\Xi_{-\eta}(x)))  
\end{array}
\end{equation}
Here, $\Delta(\Xi_{-\eta}(x),\prol u_\eps(\Xi_{-\eta}(x)))$ is negligible
by definition and Theorem \ref{factorth} shows 
$\cl[x \to Q(\eta,\Xi_{-\eta}(x),\prol u_\eps(\Xi_{-\eta}(x)))]$ 
to be moderate. It follows that $(\Delta(x,\prol(g_\eta u_\eps)(x)))_{\eps\in I}$ is
negligible, which precisely means that $g_\eta u$ is a solution in $\G$. 
\ep

\bex \label{factorex} 
Suppose that system (\ref{system}) satisfies $\frac{\partial \Delta_i}{\partial
z_{k_j}} = c_i \delta_{i k_j}$ for some $p <k_1<\dots< k_s \le N$ and 
nonzero constants $c_i$. Then (\ref{system}) satisfies conditions (i) and (ii)
of Theorem \ref{factorth}. In fact, $\det(J_{k_1,\dots,k_s}(\Delta))$ is 
constant and since the underlying domain is convex the above assumption on
(\ref{system}) is equivalent with
$$
\Delta_{i}(z) = c_i z_{k_i} + F_i(z') \qquad (1\le i\le s)
$$
with $F_i$ smooth and $z'$ as in \ref{factorth} (i). From this, injectivity of 
$z \to (z',\Delta(z))$ is immediate. Thus by \ref{factorth} any slowly
increasing classical symmetry group of (\ref{system}) is a ${\cal G}$-symmetry
group as well. This observation applies e.g.\ to the system
$$
\begin{array}{rcl}
U_t + U U_x &=& 0\\
V_t + U V_x &=& 0
\end{array}
$$
considered in \cite{KO}, Example 3.6. Also, Theorem 3.8 of \cite{KO} is a 
special case of this setup (for $s=1$).
\et

\section{Associated and distributional symmetry\\ groups}

In this section we are going to investigate symmetries of system (\ref{system}) in 
the sense of association, i.e.\ we shall be
concerned with group actions that transform
solutions of 
\begin{equation}\label{asssystem}
\Delta_i(x,u^{(n)}) \approx 0  \qquad (i = 1,\dots,s) 
\end{equation}
into solutions in the sense of association. Such group actions will be called 
{\em symmetries in the sense of association} or,
by slight abuse of terminology, {\em associated symmetries}. 
Also, we will give a rather general criterion for transferring classical symmetry
groups of linear systems to the distributional setting.
In what follows we will only consider projectable 
symmetry groups.

%%%%%%%%%%%%%%%%%%%%%%%%%%%%%%%%%%%%%%%%%%%%%%%%%%%%%%%%%%%%%%%%%%%%%%%%%%%%%
\bd
$u=(u_{1},...,u_{q}) \in
{\cal{G}}(\Omega )^{q}$ is a solution to {\rm (\ref{asssystem})} (also called 
an associated solution to {\rm (\ref{system})}) if $u$ has a representative 
$(u_{1\epsilon },...,u_{q\epsilon }) \in
{\cal{E}}_{M}(\Omega )^q$ such that for every $\varphi \in {\cal{D}}(\Omega )$
\begin{equation} \label{assdef} 
\int \Delta _{i} (x,\prol u_\eps(x))\varphi (x)\, 
dx \rightarrow 0
\mbox{ as } \epsilon \rightarrow 0, \; 1\le i \le s
\end{equation}
The set of all associated solutions to (\ref{system}) will be denoted
by ${\cal{A}} _{\Delta }$. The set of all $u \in (\G_{\infty})^q$ 
satisfying (\ref{assdef}) will be termed ${\cal{A}}{\cal{B}}_{\Delta }$.
Let ${\cal A} \subseteq {\cal A}_\Delta$.
A symmetry group $\eta\to g_\eta$ of (\ref{system}) will be called an 
${\cal{A}}$-symmetry group if $g_\eta U \in {\cal A}$ whenever $U\in {\cal A}$
and $g_\eta U$ is defined.
\et
%%%%%%%%%%%%%%%%%%%%%%%%%%%%%%%%%%%%%%%%%%%%%%%%%%%%%%%%%%%%%%%%%%%%%%

In what follows, ${\cal C}_c^k(\Om)$ $(k\in\N_0)$, the space of compactly supported
${\cal C}^k$-functions on $\Om$ is equipped with the inductive limit topology
of its subspaces ${\cal C}_K^k(\Om) = \{f\in {\cal C}^k(\Om)\mid \supp f \subseteq K\}$
($K$ compact in $\Om$).
\bd \label{strass}
Let $u = (u_{1},...,u_{q}) \in {\cal{A}} _{\Delta }(\Omega)$ 
(resp. $u \in {\cal{A}}{\cal{B}}_{\Delta }(\Omega))$ and let $k\in \N_0$. $u$
is called a $k$-strongly associated ($\stackrel{k}{\approx}$-associated)
solution to (\ref{system}) if it has a
representative $(u_{1\epsilon },...,u_{q\epsilon }) \in {\cal{E}}_{M}(\Omega )^{q}$ 
such that for every set  $B \subseteq {\cal C}^{\infty }_c (\Omega )$ which is
bounded in ${\cal{C}}_c^k(\Omega)$ we have
$$ 
\lim _{\epsilon \rightarrow 0} \sup _{\varphi \in B}
\int \Delta _{i} (x,\prol u_\eps(x))\varphi (x)\, dx = 0 \quad (1\le i\le s)
$$
The space of $k$-strongly associated solutions to {\rm (\ref{system})} is
denoted by ${\cal AS}^k_\Delta$ and we set ${\cal ABS}_\Delta^k:={\cal AS}^k_\Delta\cap
\G_\infty$. 
The corresponding symmetry groups are called ${\cal AS}^k_\Delta$-
and ${\cal ABS}^k_\Delta$-symmetry groups, respectively.
\et
%%%%%%%%%%%%%%%%%%%%%%%%%%%%%%%%%%%%%%%%%%%%%%%%%%%%%%%%%%%%%%%%%%%%%%%%%%%
\brm  
Associated solutions (i.e., solutions of type ${\cal A}_\Delta$)
play a central role in applications to numerics (cf.\ the remarks in
Section \ref{intro}). Moreover, 
by imposing slight changes in their definition,  
Colombeau type algebras
can be adapted to a wide variety of solution types. 
Thus ${\cal G}_\infty$ is particularly useful 
for the investigation of shock wave solutions of conservation laws
(cf.\ \ref{Riem}, \ref{quasiex}, \ref{3.6} below), which  
fall into this class of generalized functions. 
It will turn out in 
Theorem \ref{assth} (i) that the uniform bounddedness
of solutions is essential for the characterization of
${\cal ABS}_\Delta$-symmetry groups. 
\et
%%%%%%%%%%%%%%%%%%%%%%%%%%%%%%%%%%%%%%%%%%%%%%%%%%%%%%%%%%%%%%%%%%%%%%%%%%%
\bex \label{Riem}
Consider the Riemann problem 
\begin{equation}\label{burgerstrong}
\begin{array}{l}
u_t + uu_x = 0 \\
u(x,0) = u_l + (u_r-u_l)H(x)  
\end{array}
\end{equation}

where $H$ denotes the Heaviside function. For $u_l>u_r$ the unique weak solution
to this problem is $u(x,t)=u_l+(u_r-u_l)H(x-ct)$, $(x,t)\in \R\times [0,\infty)$,
where $c=(u_l+u_r)/2$.
Let $\theta\in \D(\R)$, $\theta\ge 0$, $\int \theta = 1$  and set 
$H_\eps(y) = \int_{-\infty}^{y/\eps} \theta(x)\,dx$, $\delta_\eps(y) =
\frac{1}{\eps}\theta(\frac{y}{\eps})$. We are looking for $1$-strongly 
associated solutions to (\ref{burgerstrong}) of the form
$u_\eps(x,t) = u_l + (u_r-u_l)H_\eps(x-ct)$, with $c$ to be determined. 
Thus let $B$ be a bounded subset of 
${\cal C}_c^1(\R\times (0,\infty))$. This means that there exists 
$K\comp \R\times (0,\infty)$
such that $\supp \varphi \subseteq K$ for all $\varphi\in B$ and
$$
\sup \{|\partial^\alpha\varphi(x,t)| \mid \vphi\in B,\,(x,t)\in K,\, |\alpha|\le 1\}
<\infty
$$
Let $\psi\in B$. Then
\begin{eqnarray*}
&& \int_{\R\times (0,\infty)} [u_{\eps t}(x,t) + u_\eps(x,t)u_{\eps x}(x,t)]\psi(x,t)
\,dx\,dt =\\  
&&  \int_{\R\times (0,\infty)} [-c(u_r-u_l)\partial_xH_\eps(x-ct) +\\ 
&&\hphantom{7em}\frac{1}{2}\partial_x((u_l+(u_r-u_l)H_\eps(x-ct))^2)] \psi(x,t)\,dxdt\\
&&
= - \int_{\R\times (0,\infty)} [-c(u_r-u_l)H_\eps(x-ct)
+\frac{1}{2}((2u_l(u_r-u_l)\\
&& \hphantom{\cdot\psi(x,t)\,dx\,dt=\dots}\cdot H_\eps(x-ct)+(u_r-u_l)^2 
H_\eps(x-ct))^2)]\psi_x(x,t)\,dx\,dt \\
&& \hphantom{\cdot\psi(x,t)\,dx\,dt=}\to \int_0^\infty [c(u_r-u_l)-
\frac{1}{2}(u_r^2-u_l^2)] \psi(ct,t)\,dt
\end{eqnarray*}
by dominated convergence, uniformly for $\psi\in B$. 
Thus $u_\eps$ as above is a $1$-strongly 
associated solution if and only if $c=(u_l+u_r)/2$, i.e.\ if and
only if the Rankine-Hugoniot jump condition is satisfied.
\et
%%%%%%%%%%%%%%%%%%%%%%%%%%%%%%%%%%%%%%%%%%%%%%%%%%%%%%%%%%%%%%%%%%%%%%%%%%%
The transfer of classical symmetry groups into symmetry groups in the sense
of association is again governed by factorization properties. 
Sufficient conditions for this transfer are provided by the following result:
%%%%%%%%%%%%%%%%%%%%%%%%%%%%%%%%%%%%%%%%%%%%%%%%%%%%%%%%%%%%%%%%%%%%%%%%%%%
\bt \label{assth}
Let $G = \{g_\eta \mid \eta \in (-\eta_0,\eta_0)\}$ 
be a slowly increasing symmetry group of {\rm (\ref{system})}
admitting a global factorization of the
form {\rm (\ref{factor0})}. Then 
\begin{itemize}
\item[(i)] If $Q$ depends exclusively on $\eta$, $x$ and $u$ then 
$g_\eta$ is an ${\cal ABS}_\Delta$-symmetry group of
(\ref{system}).
\item[(ii)] If $Q$ depends exclusively on $x$ and $\eta$
then for any $k>0$ $g_\eta$ is an ${\cal AS}_\Delta^k$-symmetry group of
(\ref{system}).
\end{itemize}
\et
%%%%%%%%%%%%%%%%%%%%%%%%%%%%%%%%%%%%%%%%%%%%%%%%%%%%%%%%%%%%%%%%%%%%%%%%%%%
\pr
(i) Let $\phi \in \D(\Xi_\eta(\Om))$ and $u \in 
{\cal ABS}_\Delta(\Om)$. It follows 
from (\ref{gafuncprgen}) that $u \in \G_{\infty}(\Xi_\eta(\Om))^q$.
Moreover, by (\ref{funcfactor}), for any  $1\le i \le s$ we have
\begin{eqnarray*}
&& \int \Delta_i(x,\prol(g_\eta u_\eps)(x))\phi(x)\,dx   \\
&& = \sum_j\int Q_{ij}(\eta,\Xi_{-\eta}(x),u_\eps(\Xi_{-\eta}(x)))
\Delta_j(\Xi_{-\eta}(x),\prol u_\eps(\Xi_{-\eta}(x)))\phi(x) \, dx\\
&& = \sum_j\int \Delta_j(x',\prol u_\eps(x'))\underbrace{Q_{ij}(\eta,x',u_\eps(x'))
\phi(\Xi_\eta(x'))|\det D\Xi_\eta(x')|}_{(*)} \, dx'\\
\end{eqnarray*}
If $\phi$ varies in a subset of $\D(\Xi_\eta(\Om))$ which is bounded in ${\cal C}_c$
then due to the fact that the $u_\eps$ are globally bounded, $(*)$ varies in
a ${\cal C}_c$-bounded subset of $\D(\Om)$. Hence the above expression tends to 
zero, uniformly in $\phi$ and $\eps$.

(ii) With $\phi$ as in (i),
\begin{eqnarray*}
&& \int \Delta_i(x,\prol(g_\eta u_\eps)(x))\phi(x)\,dx   \\
&& = \sum_j \int \Delta_j(x',\prol u_\eps(x'))
\underbrace{Q_{ij}(\eta,x')\phi(\Xi_\eta(x'))
|\det D\Xi_\eta(x')|}_{(*)} \, dx' \\
\end{eqnarray*}
which tends to zero since $(*)$ is ${\cal C}_c^k$-bounded if 
$\phi$ varies in a ${\cal C}_c^k$-bounded subset of $\D(\Xi_\eta(\Om))$.
\ep
%%%%%%%%%%%%%%%%%%%%%%%%%%%%%%%%%%%%%%%%%%%%%%%%%%%%%%%%%%%%%%%%%%%%%%%%%%%%%%
\bex \label{quasiex}
Let $A$ be an $s\times s$ matrix with ${\cal C}^1$-entries
$a^{ij}=a^{ij}(u^1,...,u^s),\; i,j$ $=$ $1,...,s$ on $\R^s$.
We consider the quasilinear system
\begin{equation} \label{quasi}
\Delta(u,u_x,u_t) \equiv u_t + A(u) u_x
\end{equation}
(for solvability resp.\ unique solvability
of (\ref{quasi}) we refer to \cite{Smol} and \cite{Bres})
and the action of a classical symmetry group
$
 G = \{
  g_\eta \mid \eta \in (-\eta_0,\eta_0)
 \}
$
\begin{equation} \label{oz}
\left.
\begin{array}{rcl}
 \tilde x &=& \Xi_\eta^1(x,t) \\
 \tilde t &=& \Xi_\eta^2(x,t) \\
 \tilde u &=& \Phi_\eta (x,t,u)
\end{array}
\right\}, \quad
 \eta \in (-\eta_0,\eta_0)
\end{equation}
where $\Xi^1_\eta,\Xi^2_\eta \in {\cal C}^\infty (\Omega),
\U_\eta \in \left({\cal C}^\infty (\Omega \times \Real^s)
\right)^s,\; \eta\in (-\eta_0,\eta_0)$. 
We are going to show that $G$ is an ${\cal AS}_\Delta^k$-symmetry group of
(\ref{quasi}) for each $k$ provided that
 \begin{enumerate}
 \item[(i)] \label{notonu} $\partial_u \Phi_\eta$ does not depend on $u$, i.e.\ 
            $G$ acts linearly on the dependent variables.
 \item[(ii)]\label{notonx} $\Xi^2_\eta$ does not depend on $x$.
 \end{enumerate}

Denoting as above the Jacobian of $\Xi_\eta$ by $J\Xi_\eta$, $\tilde u _{\tilde x}$
and $\tilde u _{\tilde t}$ are given by the first and second column of the 
matrix valued function $
  \Ui_\eta \in \left(
    {\cal C}^\infty (\Omega \times \Real^{3s})
  \right) ^{s \times 2}
$, where $\Ui_\eta (x,t,u,u_x,u_t)$ is defined by

\begin{equation}\label{xxxx}
   \left(
    \begin{bmatrix}
      \U_{\eta x} (x,t,u) &
      \U_{\eta t} (x,t,u)
    \end{bmatrix} + \frac{\partial \Phi_\eta}{\partial u} (x,t)
    \begin{bmatrix}
      u_x & u_t
    \end{bmatrix}
  \right)  (J\Xi_\eta(x,t))^{-1}
\end{equation}
Moreover,
$$
 \Delta (\tilde u,\tilde u_{\tilde x},\tilde u_{\tilde t})=:
  \tilde \Delta_\eta (x,t,u,u_x,u_t) $$
$$= \int_0^1
    (\grad{u_t}{\tilde \Delta_\eta}) (
      x,t,u,u_x,\theta u_t + (1-\theta) (-A(u)u_x)
    )
  \,d\theta \times (u_t - (-A(u)u_x))
$$
$$= Q_\eta(x,t,u) \Delta(u,u_x,u_t),
$$
Note that $Q_\eta$ does not depend on $u_x, u_t$ by our assumption on $G$
and the explicit form of $\Delta$.
Now
\begin{eqnarray}\label{Deltaeta}
\tilde\Delta_\eta =  \Ui_\eta(x,t,u,u_x,u_t) \begin{bmatrix} 0\\1 \end{bmatrix} +
  A(\U_\eta(x,t,u))
  \Ui_\eta(x,t,u,u_x,u_t) \begin{bmatrix} 1\\0 \end{bmatrix}
\end{eqnarray}
\begin{eqnarray*}
  \Ui_\eta
  \begin{bmatrix} 0\\1
  \end{bmatrix}
  &=& \left(
    \begin{bmatrix} \U_{\eta x} & \U_{\eta t}
    \end{bmatrix} + \U_{\eta u}
    \begin{bmatrix} u_x & u_t
    \end{bmatrix}
  \right) \frac 1{\det J\Xi_\eta}
  \begin{bmatrix} -\X_{\eta t} \\ \X_{\eta x}
  \end{bmatrix} \\
  &=& \frac 1{\det J\Xi_\eta} \left(
    \X_{\eta x} \U_{\eta t} - \X_{\eta t} \U_{\eta x}
    + \X_{\eta x} \U_{\eta u} u_t - \X_{\eta t} \U_{\eta u} u_x
  \right)
\end{eqnarray*}
\begin{eqnarray*}
  \Ui_\eta
  \begin{bmatrix} 1\\0
  \end{bmatrix}
  &=& \left(
    \begin{bmatrix} \U_{\eta x} & \U_{\eta t}
    \end{bmatrix} + \U_{\eta u}
    \begin{bmatrix} u_x & u_t
    \end{bmatrix}
  \right) \frac 1{\det J\Xi_\eta}
  \begin{bmatrix} \T_{\eta t} \\ -\T_{\eta x}
  \end{bmatrix} \\
  &=& \frac 1{\det J\Xi_\eta} \left(
    \T_{\eta t} \U_{\eta x} - \T_{\eta x} \U_{\eta t}
    + \T_{\eta t} \U_{\eta u} u_x - \T_{\eta x} \U_{\eta u} u_t
  \right)
\end{eqnarray*}
Thus
$$
  Q_\eta= \grad{u_t}{\tilde \Delta_\eta} = \frac 1{\det J\Xi_\eta} \left(
    \X_{\eta x} \U_{\eta u} - \T_{\eta x} A(\U_\eta) \U_{\eta u}
  \right) = \frac{
    \X_{\eta x} \eigmat_s - \T_{\eta x} A(\U_\eta)
  }{
    \X_{\eta x} \T_{\eta t} - \T_{\eta x} \X_{\eta t}
  } \U_{\eta u}
$$
where $\eigmat_s$ is the $s \times s$ identity matrix.
This implies
\begin{equation} \label{decomp1}
  \tilde \Delta_\eta (x,t,u,u_x,u_t) =
    \frac{\Xi^1_{\eta x}\eigmat_s - \Xi^2_{\eta x}A(\Phi_\eta)}
  {\Xi^1_{\eta x} \Xi^2_{\eta t}
  - \Xi^1_{\eta t} \Xi^2_{\eta x}}\Phi_{\eta u} \Delta (x,t,u,u_x,u_t) =
   Q_\eta \Delta.
\end{equation}
By assumptions (i) and (ii) it follows that $Q_\eta$ does not depend on $u$.
Thus our claim follows from \ref{assth} (ii).
\et
\bex \label{3.6}
Let $f \in {\cal C}^1 (\Real)$ and let $F$ be a primitive of $f$. The
generalized Burgers equation
\begin{equation} \label{SIGMA}
\begin{split}
& u_t + f(u) u_x = 0 \\
& u(x,0) = u_l + (u_r-u_l) H(x)
\end{split}
\end{equation}
has a weak solution of the form
\begin{equation} \label{b.1}
 u(x,t) = u_l + (u_r-u_l) H(x-ct), \quad
 (x,t) \in \Real \times (0,T)
\end{equation}
where
$
 c = (F(u_r)-F(u_l))/(u_r-u_l)
$.
(If $f'>0$ and $u_l>u_r$, then this solution is unique.)
The same arguments as in the case of Burgers' equation in
\ref{Riem} imply that
\eqref{b.1} is the $1$-strongly associated solution to \eqref{SIGMA}
if and only if $c$ is of the given form. 
In this special case of (\ref{quasi}), conditions (i) and (ii) of \ref{quasiex}
are in fact necessary for shock wave solutions to be transformed into
$1$-strongly associated solutions. In fact, if (i) or (ii) is violated then
from the explicit calculations in \ref{quasiex} it follows that the integrand 
in the proof of \ref{assth} (ii) will contain unbounded terms.

More complex factorizations arise in case $f$ is invertible. Indeed, by a 
straightforward explicit calculation, the infinitesimal generators
$$
\begin{array}{l}
{\mathbf w _1} = x t \partial _x + t^2 \partial _t +
  \frac 1{f'(u)} (x - f(u) t) \partial _u \\
{\mathbf w _2} = x^2 \partial _x + x t \partial _t +
  \frac {f(u)}{f'(u)} (x - f(u) t) \partial _u
\end{array}
$$
give rise to the group actions
\begin{itemize} \itemsep=0pt
\item[\( G_1 \):]
  \( \tilde x = x (\eta) = \frac x {1 - \eta t} \),
  \( \tilde t = t (\eta) = \frac t {1 - \eta t} \),
  \begin{math}
    \tilde u = f^{-1} (\eta x + f (u) - \eta f(u) t)
  \end{math}
\item[\( G_2 \):]
  \( \tilde x = x (\eta) = \frac x {1 - \eta x} \),
  \( \tilde t = t (\eta) = \frac t {1 - \eta x} \),
  \begin{math}
    \tilde u = f^{-1} \left(
      \frac{ f(u) }{ 1 - \eta (x - f (u) t) }
    \right)
  \end{math}
\end{itemize}
with factors 
\begin{eqnarray*}
Q_1 &=& \frac{
      (1 - \eta t) ^3 f'
    }{
      f' (f^{-1} (\eta x + f(u) - \eta f(u) t))
    }\\
Q_2 &=& \frac{
      (1 - \eta x) ^3 f' (u)
    }{
      (1 - \eta (x - t f)) ^3
      f' \left(
        f^{-1} \left(
          \frac{f (u)}{1 - \eta (x - t f(u))}
        \right)
      \right)
    }  
\end{eqnarray*}
\et
Theorem \ref{assth} raises the question of finding criteria for classical symmetry
groups to display the favorable factorization properties used above. As an important 
case where a general result is available we now turn to systems of linear PDEs. With
a view to applications in distribution theory (cf.\ section \ref{infcritsec}) we also
restrict our attention to group actions which act linearly on the dependent variables.

%%%%%%%%%%%%%%%%%%%%%%%%%%%%%%%%%%%%%%%%%%%%%%%%%%%%%%%%%%%%%%%%%%%%%%%%%%%%%%%%%%%
\bt \label{linsysfactor}
Suppose that in {\rm (\ref{system})}, $\Delta$ is a linear differential operator:
\begin{equation}\label{lindiff}
\Delta_i(z) = \sum_{k=p+1}^N a^i_k(z_1,\dots,z_p)z_k + a_0^i(z_1,\dots,z_p) 
\quad (1\le i \le s)
\end{equation}
Furthermore, let $g_\eta$ be a one parameter symmetry group of (\ref{lindiff}) which acts 
linearly in the dependent variables:
\begin{equation}\label{lingroup}
g_\eta(x,u) = (\Xi_\eta(x), \Phi_\eta(x)\cdot u + \Psi_\eta(x))
\end{equation}
If there exist $p < k_1< \dots < k_s \le N$ such that $J_{k_1,\dots,k_s}(\Delta)$ is
globally nonsingular then conditions (i) and (ii) of \ref{factorth} are satisfied
and $Q$ in {\rm (\ref{factor0})} depends exclusively on $\eta$ and
$x = (z_1,\dots,z_p)$.
\et
%%%%%%%%%%%%%%%%%%%%%%%%%%%%%%%%%%%%%%%%%%%%%%%%%%%%%%%%%%%%%%%%%%%%%%%%%%%%%%%%%%
\pr Using the same notations and conventions as in the proof of Theorem \ref{factorth}
we have $\tilde \Delta(z) = (z_1,\dots,z_{N-s},A(x)\cdot z + a_0(z))$ where the
$s\times N$ matrix $A(x)$ is of the form $(A'(x),A''(x))$ with $A'(x)$ an $s\times
(N-s)$ matrix and $A''(x) = J_{k_1,\dots,k_s}(\Delta)$ invertible.

Thus 
\begin{equation}\label{deltainv}
\tilde\Delta^{-1}(y) = (y',A''(y_1,\dots,y_p)^{-1}\cdot[y''-A'(y_1,\dots,y_p)\cdot y'
-a_0(y_1,\dots,y_p)])
\end{equation}
is affine linear in $y_k$ for $k>p$.
By (\ref{lingroup}) and \cite{Olver}, (2.18) we have
\begin{equation}\label{prolong}
\prol g_\eta(z) = (\Xi_\eta(x),\Phi_\eta(x)z_{p+1} + \Psi_\eta(x),\dots,\bar z_k,\dots,\bar z_N) \end{equation}
where
\begin{equation}\label{zk}
\bar z_k = \sum_{l=p+1}^N b^k_l(\eta,x) z_l + b_0^k(\eta,x)
\end{equation}
with certain smooth functions $b^k_l$ (Actually, the upper limit $N$ in these sums
is only attained for terms corresponding to highest order derivatives but for the
following argument only the general form of (\ref{zk}) is of interest).
Hence both $f_\eta$ and $\tilde\Delta^{-1}$ are (affine) linear in $z_k$ for $k>p$. 
It follows
that the matrix $J_{k_1,\dots,k_s}(f_\eta\circ \tilde\Delta^{-1})$ is independent of
$z_k$ for $k>p$. This observation, together with (\ref{f2}), finishes the proof. 
\ep
%%%%%%%%%%%%%%%%%%%%%%%%%%%%%%%%%%%%%%%%%%%%%%%%%%%%%%%%%%%%%%%%%%%%%%%%%%%%%%%%%%%%%%%%

\brm \label{linsysremark}
The conclusion of Theorem \ref{linsysfactor} remains valid for a semilinear
system
\begin{equation}\label{semlindiff}
\Delta_i(z) = \sum_{k=p+1}^N a^i_k(z_1,\dots,z_p)z_k + a_0^i(z_1,\dots,z_p,
z_{p+1},...,z_{p+q})
\quad (1\le i \le s)
\end{equation}
 provided that $k_1,\dots,k_s$ correspond to indices of highest order derivatives
of the dependent variables. This follows immediately from an inspection of the above proof
(only the form of $a_0$ in (\ref{deltainv}) changes).
%\end{itemize}
\et
%%%%%%%%%%%%%%%%%%%%%%%%%%%%%%%%%%%%%%%%%%%%%%%%%%%%%%%%%%%%%%%%%%%%%%%%%%%%%%%%%%%%%%%%
\bc \label{linsyscor}
Let $G$ be a symmetry group of the linear (resp.\ semilinear)
system {\rm (\ref{lindiff})} (resp.\ {\rm (\ref{semlindiff})})  such
that the assumptions of Theorem {\rm \ref{linsysfactor}} (resp.\ of Remark 
{\rm \ref{linsysremark}} ) are satisfied. Then $G$ is an ${\cal AS}_\Delta$-symmetry
group of {\rm (\ref{lindiff})} (resp. {\rm (\ref{semlindiff})}).
\et
%%%%%%%%%%%%%%%%%%%%%%%%%%%%%%%%%%%%%%%%%%%%%%%%%%%%%%%%%%%%%%%%%%%%%%%%%%%%%%%%%%%%%%%%
\pr Immediate from \ref{assth} (ii) and \ref{linsysfactor}.\ep
%%%%%%%%%%%%%%%%%%%%%%%%%%%%%%%%%%%%%%%%%%%%%%%%%%%%%%%%%%%%%%%%%%%%%%%%%%%%%%%%%%%%%%%%

\bex \label{semlin}
Let
$$
\Delta \equiv u_t+a(x,t)u_x+a_0(x,t,u)=0,\; (x,t\in \R)
$$
then setting $z=(x,t,u,u_x,u_t)$ and using the notation of \ref{linsysfactor} we have
$s=1$, $k_s=5$, $\tilde\Delta(z) = (z_1,\dots,z_4,z_5 + a(z_1,z_2)z_4 + a_0(z_1,z_2,z_3))$,
$\tilde\Delta^{-1}(y) = (y_1,\dots,y_4,y_5 - a(y_1,y_2)y_4 - a_0(y_1,y_2,y_3))$, and
$$
\begin{array}{rcl}
\mbox{pr}^{(1)}g_\eta(z)\!\! &=& \!\!(\Xi_\eta(z_1,z_2),\Phi_\eta(z_1,z_2)z_3+\Psi_\eta(z_1,z_2),
\sum_{j=3}^5 b_j^4(\eta,z_1,z_2)z_j +\\ 
&&b_0^4(\eta,z_1,z_2),
\sum_{j=3}^5 b_j^5(\eta,z_1,z_2)z_j + b_0^5(\eta,z_1,z_2))\\
\Delta(\mbox{pr}^{(1)}g_\eta(z))\!\! &=&\!\! \sum_{j=3}^5 b_j^5 z_j + b_0^5 + a\cdot(\sum_{j=3}^4 b_j^4 z_j
+ b_0^4) + a_0
\end{array}
$$
Hence
$$ Q(y) = \frac{\partial}{\partial y_5}(f_\eta\circ\tilde\Delta^{-1})(y)=  
b_5^5(\eta,y_1,y_2) + a(y_1,y_2)b_2^4(\eta,y_1,y_2).
$$
\et
%%%%%%%%%%%%%%%%%%%%%%%%%%%%%%%%%%%%%%%%%%%%%%%%%%%%%%%%%%%%%%%%%%%%%%%%%%%%%%%%%%
Turning now to the distributional setting, we first note that the most general
group actions applicable to distributions are those which are projectable and
act linearly on the dependent variables, i.e.\ those which are of the form
(\ref{lingroup}). If $u\in \D'(\Om)^q$ then the action of $g_\eta$ is defined
by 
$$
g_\eta u := \Xi_{-\eta}^*(\Phi_\eta\cdot u + \Psi_\eta)
$$
where $\Xi_{-\eta}^*$ denotes (componentwise) distributional pullback, i.e.,
$$
\langle f^*(u),\vphi\rangle = \langle u(y), \vphi(f^{-1}(y))|\det D(f^{-1})(y)|)
\rangle \qquad (u\in \D'(\Om'), \vphi \in \D(\Om))
$$
for $f: \Om \to \Om'$ a diffeomorphism. We will sometimes also write $u\circ f$ 
instead of $f^*u$.

\bd \label{distsymmdef}
Suppose that (\ref{system}) is linear and let $G$ be a local transformation group
acting linearly on the dependent variables. $G$ is called a distributional (or $\D'$-)
symmetry group of (\ref{system}) if it transforms distributional solutions of 
(\ref{system}) into distributional solutions.
\et
%%%%%%%%%%%%%%%%%%%%%%%%%%%%%%%%%%%%%%%%%%%%%%%%%%%%%%%%%%%%%%%%%%%%%%%%%%%%%%%%%%%%%%%%
\ref{distsymmdef} is the most general definition of distributional symmetry groups.
More restrictive notions (as introduced e.g.\ in \cite{ga}) will be discussed in
the following section.
%%%%%%%%%%%%%%%%%%%%%%%%%%%%%%%%%%%%%%%%%%%%%%%%%%%%%%%%%%%%%%%%%%%%%%%%%%%%%%%%%%%%%%%%
\bt \label{distsymmth}
Let $G$ be a symmetry group of the linear system {\rm (\ref{lindiff})} such 
that the assumptions of Theorem {\rm \ref{linsysfactor}}  
are satisfied. Then $G$ is a distributional symmetry group of {\rm (\ref{lindiff})}.
\et
%%%%%%%%%%%%%%%%%%%%%%%%%%%%%%%%%%%%%%%%%%%%%%%%%%%%%%%%%%%%%%%%%%%%%%%%%%%%%%%%%%%%%%%%
\pr By \ref{linsysfactor} and (\ref{funcfactor}), for any smooth 
$u:\Om \to \R^q$ we have 
$$
%\begin{array}{l}
\Delta(x,\prol(g_\eta u)(x))
=Q(\eta,\Xi_{-\eta}(x)) \cdot 
\Delta(\Xi_{-\eta}(x),\prol u(\Xi_{-\eta}(x)))  
%\end{array}
$$
Now suppose that $u\in \D'(\Om)^q$ is a solution to (\ref{system}). Choose some 
$\phi \in \D(\R^p)$ with $\int \phi(x) \, dx = 1$ and set $u_\eps^i = u_i*\phi_\eps$,
$u_\eps = (u_\eps^1,\dots,u_\eps^q)$. Then $u_\eps^i$ is smooth and converges to 
$u_i$ in $\D'$ for $\eps\to 0$. Since $g_\eta$ acts linearly on the dependent variables 
we also have that $\prol (g_\eta u_\eps) \to \prol (g_\eta u)$ in $\D'$.
Hence
\begin{eqnarray*}
&& \Delta(x,\prol(g_\eta u)(x)) =
 \lim_{\eps\to 0} Q(\eta,\Xi_{-\eta}(x)) \cdot \Delta(\Xi_{-\eta}(x),
\prol u_\eps(\Xi_{-\eta}(x))) \\
&& \hphantom{\Delta(x,\prol(g_\eta u)(x))}
= Q(\eta,\Xi_{-\eta}(x)) \cdot \Delta(\Xi_{-\eta}(x),\prol u(\Xi_{-\eta}(x)))
= 0
\end{eqnarray*}
in $\D'$, which concludes the proof.
\ep

\section{Infinitesimal criteria} \label{infcritsec}
In this section we develop infinitesimal criteria for finding symmetries 
applicable to all three settings of interest (distributional, weak, Colombeau).
The results introduced here also establish a direct connection
of our approach to symmetries of distributional and weak solutions to 
that given in \cite{ga} (cf.\ also \cite{huy,wi,bi}). 

In \cite{ga}, systems of the form
\begin{equation}\label{berlin}
 \begin{bmatrix}
  L^{11} &
  \ldots &
  L^{1q} \\
  \hdotsfor 3 \\
  L^{s1} &
  \ldots &
  L^{sq}
 \end{bmatrix}
 \begin{bmatrix}
  u_1 \\
  \hdotsfor 1 \\
  u_q
 \end{bmatrix} =
 \begin{bmatrix}
  F_1 \\
  \hdotsfor 1 \\
  F_s 
 \end{bmatrix}, 
 \text{or, for short}
\quad
 L(x,D) u = F \quad
 \end{equation}
where
$$
 L^{ij} (x,D) =
 \sum_{\abs|J|=0}^{n} a_J^{ij} (x) D^J, \quad
 a_J^{ij} \text{ smooth,} \quad
 i=1,\ldots,s, \;
 j=1,\ldots,q
$$
($n$ the order of $L$) and $F= x\mapsto F(x)\in \D'(\Om)$, $\Om\subseteq \R^p$, 
are examined. The form (\ref{berlin}) of writing a system of linear PDEs 
(which, of course, is equivalent to (\ref{lindiff})) provides the advantage of
allowing to derive very concise forms of infinitesimal criteria for 
factorization properties, cf.\ (\ref{berestinf}) below.

In Berest's approach, symmetry groups of (\ref{berlin}) are  {\em defined} via 
factorization properties: Let $G = \{g_\eta \mid \eta\in
(-\eta_0,\eta_0)\}$ be a projectable one parameter group acting linearly on
the dependent variables. Thus  $g_\eta(x,u) = (\Xi(\eta,x),\Phi(\eta,x,u))$ 
is of the form  ($i=1,\dots,p$, $k,l=1,\ldots,q$) 
$$
 \tilde x^i = \Xi^i(\eta,x),\
 \tilde u^k = \Phi(\eta,x,u) =
 \sum_{l=1}^q \varphi^{kl} (\eta,x) u^l + \psi^k (\eta,x)
$$
$G$ is called a symmetry group of (\ref{berlin}) if there exists a smooth
matrix valued map $(\eta,x) \to Q(\eta,x)$ such that (see \cite{ga}, (1.6)):
\begin{equation} \label{star}
L (\tilde x, D) (g_\eta u)(\tilde x) - F(\tilde x) =
  Q (\eta,x) (L(x,D)u(x) - F(x))\,.
\end{equation}
for all $u\in \D'(\R^p)^q$
Note that this form corresponds exactly to (\ref{funcfactor}) with $Q$ depending on 
$\eta$ and $x$ exclusively.

\cite{ga}, (1.7) gives the following infinitesimal criterion for the validity of (\ref{star}):
\begin{equation} \label{berestinf}
 [\xi D,L]u + L (
  \alpha (x) u + \beta (x)
 ) - \xi D F =
 \frac{\partial Q}{\partial \eta} \Big|_{\eta=0} (Lu-F)
\end{equation}
where
$$
 \xi = \frac{\partial \Xi}{\partial \eta} \Big|_{\eta=0},
 \quad \xi D = \sum_{i=1}^p \xi_i\frac{\partial}{\partial x_i}\,,
 \quad
 \alpha = \frac{\partial \varphi}{\partial \eta} \Big|_{\eta=0}
 \,, \quad
 \beta = \frac{\partial \psi}{\partial \eta} \Big|_{\eta=0}
$$
The equivalence of (\ref{star}) and (\ref{berestinf}) raises the question whether
factorization properties for general systems (\ref{system}) of differential equations 
can always be characterized by infinitesimal conditions similar to (\ref{berestinf}). 
For $\Delta$ smooth, an affirmative answer is given by the following result,
contained implicitly in \cite{Olver} (cf.\ \ref{factorinfrm} (ii) below). 
Since its method of proof will form the basis for our generalizations to the 
$\G$- resp.\ $\D'$-settings we state it explicitly.
\bp \label{factorinfprop}
Let $G=\{g_\eta \mid \eta\in (-\eta_0,\eta_0)\}$ be a one parameter group 
with infinitesimal generator $\mathbf{v}$ acting on ${\cal M}$. 
Then the following are equivalent
\begin{itemize}
\item[(i)] There exists a smooth mapping $Q: {\cal W}
\to \R^{s^2}$ such that 
\begin{equation}\label{factorinf1}
\Delta(\prol g_\eta(z)) = Q(\eta,z)\cdot \Delta(z) \quad 
((\eta,z)\in {\cal W})
\end{equation}
\item[(ii)] There exists a smooth mapping $\tilde Q: {\cal M}^{(n)}
\to \R^{s^2}$ such that 
\begin{equation}\label{factorinf2}
\prol \mathbf{v}(\Delta)(z) = \tilde Q(z)\cdot \Delta(z) \quad
(z\in {\cal M}^{(n)})
\end{equation}
\end{itemize}
\et
\pr (i) $\Rightarrow$ (ii): Noting that $\prol \mathbf{v}$ ($\mathbf v$ 
the infinitesimal generator
of $\eta\to g_\eta$) is a vectorfield on ${\cal M}^{(n)}$ whose flow is precisely
$\prol g_\eta$, differentiation of (\ref{factorinf1}) with respect to 
$\eta$ at $\eta=0$
gives (\ref{factorinf2}) (with $\tilde Q(z) = \frac{\partial }{\partial \eta} 
\big|_{\eta=0} Q(\eta,z)$).

(ii) $\Rightarrow$ (i): 
(\ref{factorinf2}) yields
the following linear ODE for $\Delta\circ\prol g_\eta$:
\begin{equation} \label{ODE}
\begin{array}{rcl}
\frac{\partial}{\partial \eta}(\Delta\circ \prol g_\eta(z)) &=& 
\tilde Q(\prol g_\eta(z))\cdot \Delta\circ\prol g_\eta(z)\\
\Delta\circ \prol g_\eta(z)|_{\eta=0} &=& \Delta(z)
\end{array}
\end{equation}
Let $Q(\eta,z)$ be a principal matrix solution to (\ref{ODE})  
(i.e.\  the $i$-th column of $Q$ is precisely the solution with initial value
$e_i$). Then we immediately obtain the unique solution to (\ref{ODE})
in the form
$$
\Delta\circ \prol g_\eta(z) = 
Q(\eta,z)\cdot \Delta(z)
$$
\ep
\brms \label{factorinfrm} 
\begin{itemize}
\item[(i)] From the explicit formulae given in the above proof it follows that 
$Q$ depends exclusively on $z\in {\cal M}^{(k)}$ for $k<n$ (e.g. exclusively
on $(z_1,\dots,z_p) = x$) if and only if the same is true of $\tilde Q$: 
indeed it suffices 
to note that by \cite{Olver}, (2.20) $\pi^n_k\circ\prol g_\eta = 
\mbox{\rm pr}^{(k)} g_\eta$ (where $\pi^n_k: {\cal M}^{(n)} \to {\cal M}^{(k)}$
is the natural projection).
\item[(ii)] By \cite{Olver}, eq.\ (2.26), for any nondegenerate system
(\ref{system}), (\ref{factorinf2}) is equivalent
to $\prol \mathbf{v}(\Delta_i) = 0$ on the zero-set of $\Delta$ ($1\le i \le s$)
which in turn (by \cite{Olver}, Th.\ 2.71) is necessary and sufficient for 
$\mathbf{v}$ to generate a one-parameter symmetry group of (\ref{system}). Thus
(\ref{factorinf2}) is precisely the infinitesimal version of the global 
factorization (\ref{factor0}).
\item[(iii)] It follows immediately from the definition of $\prol g_\eta$ 
that (\ref{factorinf1}) is equivalent with
\begin{equation}\label{factorinf1'}
  \begin{array}{l}
\Delta(\Xi_\eta(x,u(x)),\prol g_\eta(u)(\Xi_\eta(x,u(x)))) \\
= Q(\eta,x,\prol u (x))\cdot \Delta(x,\prol u(x)) \quad 
\forall u \in {\cal C}^\infty(\Om)^q \ \forall x\in \Om\,
  \end{array}
\end{equation}
where $\Om\subseteq \R^p$ runs through all open sets. A similar reformulation
is valid for  (\ref{factorinf2}).
\end{itemize}
\et
In the Colombeau setting the general form of (\ref{system}) allows for $\Delta$ 
itself to be a generalized function. More precisely, we suppose that 
$\Delta \in \G_\tau({\cal M}^{(n)})^s$. Thereby, the admissible symmetry
transformations themselves will become generalized functions, so-called 
(projectable)
{\em generalized group actions} $g \in \tilde \G_\tau(\R\times {\cal M})^{p+q}$.
Thus $g$ is supposed to satisfy $\G_\tau$-bounds with respect to the group
parameter $\eta \in \R$ and $\G$-bounds with respect to $(x,u)\in {\cal M}$
(\cite{KO}, Def.\ 4.8). 
Many of the infinitesimal methods of classical group analysis can be recovered in
this setting. For a detailed analysis we refer to \cite{KO}. The
proof of the analogue to \ref{factorinfprop} in the present situation 
requires the following auxiliary result on solutions of linear ODEs in $\G$:

\blem \label{godelem}
Let $A\in \tilde\G_\tau((-\eta_0,\eta_0)\times\R^p)^{m^2}$ 
such that for all $K\subset\subset \R^p$,
$$\sup_{x\in K} \int_\R \|A_\eps(s,x)\|\,ds = O(|\log(\eps)|)$$ 
Then for each $u_0\in \G(\R^p)^m$ the initial value problem
\begin{equation}\label{gode}
\begin{array}{rcl}
\partial_t u(t,x) &=& A(t,x)u(t,x)\\
u(0,x) &=& u_0(x)
\end{array}
\end{equation}
has a unique solution $u$ in $\tilde\G_\tau((-\eta_0,\eta_0) \times \R^p)^m$. 
Setting $U$ the matrix with columns
the unique solutions with initial conditions $e_i$ ($1\le i\le m$), the solution
to (\ref{gode}) is given by $U(t,x)u_0(x)$. We call $U$ a principal matrix 
solution to (\ref{gode}).
\et
\pr We only sketch the argument (for details, cf.\ \cite{HO, Ligeza1, Ligeza2}). 
Choosing representatives
$(A_\eps)_\eps$ of $A$ and $(u_{0\eps})_\eps$ of $u_0$, by the corresponding
result in the ${\cal C}^\infty$-setting we obtain representatives $(u_\eps)_\eps$,
$(U_\eps)_\eps$ satisfying the claimed properties for each fixed $\eps$.
Then
$$
\|u_\eps(t,x)\| \le \|u_\eps(0,x)\| + \int_0^t \|A_\eps(s,x)\| \|u_\eps(s,x)\|\,ds
\,, \ t\in (-\eta_0,\eta_0)
$$
Thus Gronwall's inequality and the supposed growth restriction on $A$ 
imply moderateness of the representatives and unique solvability.
\ep\ms
Note that the corresponding statement for initial value
problems in the sense of association (i.e.\ replacing $=$ by $\approx$ 
in (\ref{gode})) is false since
unique solvability of linear ODEs breaks down in that context. As an easy example
take $u_\eps(x) = \frac{1}{\eps}\cos(\eps^2x)$ ($x\in \R$). Then $u' \approx 0$
but $u$ is not associated to any constant.

An element $u$ of $\G(\Om)$ is called {\em globally of $L^\infty$-log-type}
(cf.\ \cite{HO}) if 
it possesses a representative $(u_\eps)_\eps$ with 
$\sup_{x\in\Omega}|u_\eps(x)| = O(|\log(\eps)|)$.
After these preparations we can state (for the notion of 
$\G$-n-completeness, see \cite{KO}, Def.\ 4.15):

\bp \label{factorinfpropG}
Let $\Delta \in \G_\tau({\cal M}^{(n)})^s$ and let $g$ be a 
$\G$-n-complete group action on ${\cal M}$.
Consider the statements
\begin{itemize}
\item[(i)] There exists $Q \in (\tilde \G_\tau({\cal W}))^{s^2}$  
with
\begin{equation}\label{factorinfG1}
\Delta(\prol g_\eta(z)) = Q(\eta,z)\cdot \Delta(z) \ \mbox{ in } 
\tilde \G_\tau({\cal W})^{s}
\end{equation}
\item[(ii)] There exists $\tilde Q \in  {\cal G}({\cal M}^{(n)})^{s^2}$
\begin{equation}\label{factorinfG2}
\prol \mathbf{v}(\Delta)(z) = \tilde Q(z)\cdot \Delta(z) \ \mbox{ in }
{\cal G}({\cal M}^{(n)})^s
\end{equation}
\end{itemize}
Then (i) implies (ii).
If $\tilde Q\circ \prol g_\eta$ is an element of $\tilde
\G_\tau({\cal W})^{s^2}$ 
satisfying the growth property given in \ref{godelem} then (ii) implies (i).
If $Q$ and $\tilde Q$ are supposed to be globally of $L^\infty$-log-type then
(i) and (ii) are equivalent.
\et
\pr Using \ref{godelem}, the proof proceeds along the lines of 
\ref{factorinfprop}.\ep\ms
From the pointvalue characterization of Colombeau generalized functions given
in \cite{KO} it follows that (\ref{factorinfG1}) is equivalent with
\begin{equation}\label{gfactorfunc}
  \begin{array}{l}
\Delta(\Xi_\eta(x),\prol g_\eta(u)(\Xi_\eta(x))) \\
= Q(\eta,x,\prol u (x))\cdot \Delta(x,\prol u(x)) \quad 
\mbox{ in } \G(\Om)^s \ \forall u \in {\cal G}(\Om)^q \ \forall \Om
  \end{array}
\end{equation}
(cf.\ \cite{KO}, Lemma 4.13 and Prop.\ 4.14).

\brm
If in \ref{godelem}, $A\in {\cal C}^\infty((-\eta_0,\eta_0);{\cal O}_C(\R^n))^{m^2}$
and $u_0 \in {\cal O}_C(\R^n)$,
then the solution to (\ref{gode}), $U(t,x)u_0(x)$ belongs to 
${\cal C}^\infty((-\eta_0,\eta_0);{\cal O}_C(\R^n))^{m}$. Thus, if we suppose
$\Delta\in {\cal O}_C(\R^N)^s$ in \ref{factorinfpropG} then we have 
$Q\in {\cal C}^\infty((-\eta_0,\eta_0);{\cal O}_C(\R^N))^{s^2}$,
$\tilde Q \in {\cal O}_C(\R^N)^{s^2}$ and (i) and (ii) in \ref{factorinfpropG}
are equivalent
\et
Turning now to the distributional setting, we first note the following result on
linear ODEs in $\D'$.

\blem \label{distodelem}
Let $A\in {\cal C}^\infty(\R^p)^{m^2}$ and let $a\in \D'(\R^p)^m$.
Then the initial value problem 
\begin{equation}\label{distode}
\begin{array}{rcll}
\partial_t u &=& A\cdot u \qquad & \mbox{in} \ {\cal C}^\infty(\R,\D'(\R^p)^m)\\
u(0,.) &=& a \qquad & \mbox{in} \ \D'(\R^p)^m
\end{array}
\end{equation}
has the unique solution $Q\cdot a$ where $Q$ is the smooth principal matrix solution
of (\ref{distode}).
\et
\pr
That $Q\cdot a$ is a solution follows easily by regularizing the
initial data via convolution with a standard mollifier and then
taking the distributional limit of the solutions to the resulting smooth
problems. Uniqueness follows from uniqueness of the 
corresponding smooth
initial value problem by observing that for any solution $u\in 
{\cal C}^\infty(\R,\D'(\R^p)^m)$ of (\ref{distode}) 
with $a=0$ we have $\partial_t(Q^{-1}u) = 0$. But then $u = Q\cdot c$ 
with $c$ a constant vector which is necessarily $0$. 
\ep\ms
Let us suppose that 
$\Delta$ is of the
form (\ref{lindiff}) with $a^i_0\in \D'(\R^p)$ 
and $a^i_k \in {\cal C}^\infty(\R^p)$ ($1\le i \le s$, $p+1\le k\le N$). 
This is the most general 
form of differential operators applicable to elements $u$ of $\D'(\R^p)$.
Moreover, we suppose that the group action $g_\eta$ is of the form (\ref{lingroup})
(also the most general action applicable to distributions).
We consider $\Delta$ as an element of $\D'(\R^N)$ by embedding $a_0$ as
$a_0\otimes 1_{N-p}$ into $\D'(\R^N)$. 
For the following result, to simplify notations we suppose that ${\cal M}=
\R^p\times \R^q$ and that $g$ is defined on all of ${\cal M}$. 
\bp \label{fipdistr}
Under the assumptions formulated before \ref{distodelem}, 
the following are equivalent
\begin{itemize}
\item[(i)] There exists a smooth mapping $Q:\R\times \R^p \to \R^{s^2}$ 
such that 
\begin{equation}\label{factorinfdist1}
(\prol g_\eta)^*\Delta = Q(\eta,\,.\,)\cdot \Delta \ 
\mbox{ in } {\cal C}^\infty(\R,\D'({\cal M}^{(n)})^s) 
\end{equation}
\item[(ii)] There exists a smooth mapping $\tilde Q: \R^p
\to \R^{s^2}$ such that 
\begin{equation}\label{factorinfdist2}
\prol \mathbf{v}(\Delta) = \tilde Q \cdot \Delta \ \mbox{ in } \D'({\cal M}^{(n)})^s
\end{equation}
\item[(iii)] There exists a smooth mapping $Q:\R\times \R^p \to \R^{s^2}$ such that
\begin{equation}\label{fd}
  \Delta(\Xi_\eta(.),\prol(g_\eta u)) = Q(\eta,\,.\,)\cdot\Delta(.,\prol u(.)) 
\ \mbox{ in } 
{\cal C}^\infty(\R,\D'(\R^p)^s) 
\end{equation}
for all $u\in \D'(\R^p)$.
\end{itemize}
\et
\pr (i) $\Leftrightarrow$ (ii): 
Using \ref{distodelem}, the proof is again identical to that of 
\ref{factorinfprop} (for the distributional identity
$\frac{\partial}{\partial \eta}((\prol g_\eta)^*\Delta) = \prol \mathbf{v}(\Delta)$
used in the argument, see e.g., \cite{marsden}, Th.\ 3.7).

(i) $\Leftrightarrow$ (iii): With the notations introduced in (\ref{lindiff}),
(\ref{prolong}), (\ref{zk}) and $z=(x_1,\dots,x_p,$ $z_{p+1},\dots,z_N)$ we have
\begin{eqnarray*}
\Delta_i(\prol g_\eta(z)) &=& \sum_{j=p+1}^Na_j^i(\Xi(\eta,x))\left(\sum_{l=p+1}^N
b_l^j(\eta,x)z_l + b_0^j(\eta,x)\right) + \\
&& a_0^i(\Xi(\eta,x))\otimes 1_{N-p} \quad (i=1,\dots,s)
\end{eqnarray*}
It follows that (\ref{factorinfdist1}) can be written in the form
\begin{eqnarray*}
&&\sum_{l=p+1}^N \left(\sum_{j=p+1}^N a^i_j(\Xi(\eta,x))b^j_l(\eta,x)\right)z_l +
\sum_{j=p+1}^N a^i_j(\Xi(\eta,x))b_0^j(\eta,x) + \\
&& a_0^i(\Xi(\eta,x))\otimes
1_{N-p} = \sum_{m=1}^s Q_{im}(x)\left(\sum_{k=p+1}^N a^m_k(x)z_k + a_0^m(x)
\otimes 1_{N-p}\right)  
\end{eqnarray*}
Hence, introducing suitable smooth functions $\tilde a^i_k$ and distributions
$\tilde a_0^i$ ($k=p+1,\dots,N$, $i=1,\dots,s$), the proof reduces to establishing
the equivalence of
\begin{eqnarray*}
&& \langle \sum_{k=p+1}^N \tilde a^i_k(\eta,x) z_k + \tilde a_0^i(x)\otimes 1_{N-p},
\vphi_1(x_1)\dots\vphi_p(x_p)\psi_{p+1}(z_{p+1})\dots\psi_N(z_N)\rangle\\ 
&& =0 \,, \quad
\vphi_m,\psi_n\in \D(\R)\,, \, m=1,\dots,p,\, n=p+1,\dots,N
\end{eqnarray*}
and
\begin{eqnarray*}
&& \langle \sum_{k=p+1}^N \tilde a^i_k(\eta,x) z_k(x) + \tilde a_0^i(x)\otimes 1_{N-p},
\vphi_1(x_1)\dots\vphi_p(x_p)\rangle = 0\,,\\ 
&& 
\vphi_m\in \D(\R)\,, \, m=1,\dots,p,\, z_k\in {\cal C}^\infty(\R^p),\, k=p+1,\dots,N
\end{eqnarray*}
This last assertion naturally splits into a (distributional) invariance property of
$\tilde a_0^i\otimes 1_{N-p}$ resp.\ $\tilde a_0^i$ and a smooth part, 
both of which are easily seen to be equivalent.
\ep\ms
Thus (\ref{star}) corresponds precisely to (\ref{fd}). In particular, the specific 
form of
the infinitesimal criterion for the validity of (\ref{star}) follows from an explicit
calculation of $\prol \mathbf{v}(\Delta)$ in the notation (due to Berest) introduced
at the beginning of this section. In fact, we have
$$
\mathbf{v} = \sum_{i=1}^p \xi^i(x) \partial_{x_i} + \sum_{k=1}^q
\underbrace{\left(\sum_{l=1}^q \alpha^{kl}(x)u^l + \beta^k(x)\right)}_{=:\Phi_k(x,u)}\partial_{u_k}
$$  
Then by \cite{Olver}, Th.\ 2.36, $\prol \mathbf{v} = \sum_{i=1}^p\xi^i\partial_{x_i} 
 +  \sum_{k=1}^q
\sum_{|J|=0}^n \Phi_k^J(x,u^{(n)}) \partial_{u^k_J}$, where 
$$
\Phi^J_k(x,u^{(n)}) = D_J\left(\Phi_k - \sum_{i=1}^p \xi_i\frac{\partial u^k}{\partial x^i}
\right) + \sum_{i=1}^p \xi_i u^k_{J,i}
$$
($u^k_{J,i} = \frac{\partial u^k_j}{\partial x_i}$). Thus the $r$-th
component of $\prol \mathbf{v}(\Delta) = \prol \mathbf{v}(L(x,D)u-F)$ is calculated 
as follows:
\begin{eqnarray*}
&& \prol \mathbf{v}\left(\sum_{k=1}^q\sum_{|J|=0}^n a_J^{rk} u^k_J\right) - 
\prol \mathbf{v}(F_r) \\
&& =\sum_{k=1}^q\sum_{|J|=0}^n\sum_{i=1}^p \xi_i\partial_i(a^{rk}_J) u^k_J 
+ \sum_{k=1}^q\sum_{|J|=0}^n \Phi^J_k a^{rk}_J  
- \sum_{i=1}^p \xi_i \partial_i F_r \\
&& = \sum_{k=1}^q\sum_{|J|=0}^n\sum_{i=1}^p \xi_i\partial_i(a^{rk}_J) u^k_J
- \sum_{k=1}^q\sum_{|J|=0}^n\sum_{i=1}^p a^{rk}_J\partial^J\xi_i \partial_i u^k  \\
&& + \sum_{k=1}^q\sum_{|J|=0}^n\sum_{l=1}^q a^{rk}_J \partial^J \alpha^{kl} u^l 
+ \sum_{k=1}^q\sum_{|J|=0}^n\sum_{l=1}^q \alpha^{kl}  a^{rk}_J\partial^J u^l \\
&& + \sum_{k=1}^q\sum_{|J|=0}^n a^{rk}_J \partial^J\beta^k  
- \sum_{i=1}^p \xi_i \partial_i F_r \\
&& = \left([\xi D,L]u + L(x,D)(\alpha u + \beta) - \xi D F\right)_r
\end{eqnarray*}
This calculation, combined with \ref{fipdistr} and (\ref{fd}) provides a
rigorous proof of the equivalence of (\ref{star}) and (\ref{berestinf}).
Moreover, \ref{factorinfprop} allows
to derive infinitesimal criteria for factorization properties even for  systems that are
not necessarily linear and thereby to obtain workable criteria for finding symmetries
of weak or Colombeau solutions of such systems. For example, let us consider the
semilinear system
\begin{equation}\label{semilinear}
 L(x,D)u=F(u)
\end{equation}
where we shall suppose $F\in {\cal O}_M(\R^q)$ (to allow for an insertion of Colombeau
functions, note however that \ref{semlininfprop} below does not use this assumption). 
Furthermore, let us assume that the group action $g_\eta$ is of the form
(\ref{lingroup}). Then we have
\bp \label{semlininfprop} Under the above assumptions, the following are equivalent:
\begin{itemize}
\item[(i)] There exists a smooth mapping $Q:{\cal W}
\to \R^{s^2}$ such that for all $u\in {\cal C}^\infty(\Om)$ $(\Om\subseteq \R^p)$
and all $x\in \Om$
\begin{equation}\label{slf1}
L(\tilde x,D)(g_\eta u)(\tilde x) - F(g_\eta u(\tilde x)) = Q(\eta,x)(L(x,D)u(x)+F(u(x)))
\end{equation}
\item[(ii)] There exists a smooth mapping $\tilde Q: {\cal M}^{(n)}
\to \R^{s^2}$ such that for all $u\in {\cal C}^\infty(\Om)$ $(\Om\subseteq \R^p)$
and all $x\in \Om$
\begin{equation}\label{semlininf}
\begin{array}{l}
[\xi D,L]u(x) + L(x,D)(\alpha(x)u(x) + \beta(x)) - J(F)(u(x))\cdot \\
(\alpha(x)\cdot u(x) + \beta(x)) = \tilde Q(x) (L(x,D)u(x)-F(u(x)))
\end{array}
\end{equation}
(with $J(F)$ the Jacobian of $F$).
\end{itemize}
\et
\pr By \ref{factorinfprop} and \ref{factorinfrm} (iii) it suffices to calculate 
$\prol \mathbf{v}(\Delta)$ for $\Delta = L(x,D)u-F(u)$ as above. Noting that
$\prol \mathbf{v}(F) = J(F)(u)\cdot (\alpha\cdot u + \beta)$, the result follows
exactly as in the above calculation.\ep\ms
In \cite{ga}, it was shown that a certain splitting of (\ref{berestinf}) is
advantageous for establishing a connection between determining symmetry groups
of PDEs and group invariance of the solutions themselves, in particular with
a view to determining group invariant fundamental solutions of linear 
systems. In our more general setup, we first note that (using obvious abbreviations)
setting $K = \xi D$ $-$ $\alpha$ and $h = \alpha - 
\frac{\partial}{\partial\eta}|_{\eta=0}$,
(\ref{berestinf})  is  equivalent to (\ref{berestinf1}) as well as to
(\ref{berestinf2}). Also, (\ref{semlininf}) is equivalent to
(\ref{semlininf1}):

\begin{eqnarray}
&& (-[L,K] + hL)u - (K+h)F + L\beta = 0 \label{berestinf1} \\
&& [L,K] = hL \ , \ (K+H)F  = L\beta \label{berestinf2} \\
&& (-[L,K] + hL)u - J(F)(u)\alpha u + 
\frac{\partial}{\partial\eta}\Big|_{\eta=0}Q F(u) 
-J(F)(u)\beta = 0 \label{semlininf1} 
\end{eqnarray}
We first note the
following immediate consequence of (\ref{berestinf}):

\bp
Let $\eta \to g_\eta$ be a one-parameter symmetry group of the form (\ref{lingroup}) 
of the homogeneous system $L(x,D)u = 0$ satisfying (\ref{star}). 
If $u\in {\cal G}(\Om)^q$ is a solution,
resp.\ associated solution resp.\ strongly associated solution then so is $Ku$.
\ep
\et
\bd
Let $\eta \to g_\eta$ be a slowly increasing one-parameter group. 
We say that $u\in {\cal G}(\Om)^q$ is ${\cal G}$-invariant under $g_\eta$, 
$\approx$-invariant,
or $\stackrel{k}{\approx}$-invariant, respectively, if $g_\eta u = u$, 
$g_\eta u \approx u$, or $g_\eta u \stackrel{k}{\approx} u$ for all $\eta$.
If $g_\eta$ is of the form (\ref{lingroup}) and $u\in \D'(\Om)^q$ then $u$
is called $\D'$-invariant under $g_\eta$ if $g_\eta u = u$ in $\D'$.
\et

\bp
Let $\eta\to g_\eta$ be of the form (\ref{lingroup}). 
\begin{itemize}
\item[(i)]
Let $u \in ({\cal C}^\infty)^q$ resp.\ 
$u\in {\cal G}^q$ resp.\ $u\in (\D')^q$. 
A necessary and sufficient condition for $u$ to 
be invariant resp.\ ${\cal G}$-invariant 
resp.\ $\D'$-invariant under $g_\eta$ 
is that $Ku$ equals
$\beta$ in ${\cal C}^\infty$ resp.\  
in 
${\cal G}$ 
resp.\ in $\D'$.  
\item[(ii)] 
Then $u$ is $\approx$-invariant,
resp.\ $\stackrel{k}{\approx}$-invariant,
if $Ku - \beta \approx 0,$
resp.\ $Ku - \beta \stackrel{k}{\approx} 0$ ($\;k\in \N_0$).
\end{itemize}
\et
\pr (i) We first note that in each of the possible settings we have
$$
\frac{\partial}{\partial \eta}\Big|_0(g_\eta u) = -Ku + \beta
$$
which is immediate from the chain rule. Thus the conditions are necessary.
Conversely, for fixed $u$ set
$$
f(\eta,x) = \Xi_\eta^* u(x) - \phi(\eta,x)u(x) - \psi(\eta,x)
$$
Then $f'(\eta,x) = \Xi_\eta^*\alpha \cdot f(\eta,x)$ and $f(0,x)=0$.
Thus the claim follows from unique solvability of linear ODEs in each of the 
respective settings (cf.\ \ref{godelem} with $A$ smooth (hence automatically
satisfying the necessary growth restrictions) and \ref{distodelem}).

(ii) 
Setting $r_\eps(\eta,\,.\,) = \Xi_\eta^*(Ku_\eps - \beta)$
we have 
\begin{eqnarray*}
f_\eps'(\eta,x) &=& \Xi_\eta^*\alpha \cdot f_\eps(\eta,x) +r_\eps(\eta,x)\\
f_\eps(0,x) &=& 0
\end{eqnarray*}
Also, $r_\eps(\eta,\,.\,)\approx 0$ resp.\ $\stackrel{k}{\approx} 0$ for all
$\eta$ follows from our assumption (by substituting for $\Xi_\eta(x)$ in the
respective integrals). 
With $U(\eta,x)$ the principal matrix solution to the corresponding homogeneous
system we obtain the solution to this initial value problem in the form
$$
f_\eps(\eta,x) = U(\eta,x)\int_0^\eta U(\eta',x)^{-1}r_\eps(\eta,x)\,d\eta'
$$
Thus for $\varphi\in \D,$ resp.\ $\varphi \in B\subseteq {\cal C}^\infty_c,\;
B$ bounded in ${\cal C}^k_c$, we have
$$
\int f_\eps(\eta,x) \varphi(x)\,dx = 
\int r_\eps(\eta,x) \left(U(\eta,x)\int_0^\eta U(\eta',x)^{-1}\,d\eta'
\varphi(x)\right)\,dx
\to 0  \\
$$
as $\eps\to 0$, resp.\ this limit is uniform for  $\varphi \in B$.
It follows that $f_\eps(\eta,.)\approx 0$,
resp. $f_\eps(\eta,.) \stackrel{k}{\approx}0$
i.e., that $g_\eta u_\eps \approx u_\eps,$ resp.
$g_\eta u_\eps\stackrel{k}{\approx}u_\eps.$
\ep

Note that the converse assertion in (ii) does not follow since
$g_\eta u -u \approx 0$ resp.\ $\stackrel{k}{\approx} 0$ for all $\eta$ does 
not imply any information on $\frac{d}{d\eta}\big|_0 g_\eta u = 
-Ku_\eps + \beta$.

\bex \label{infex1}
We derive infinitesimal criteria for symmetries (\ref{oz}) of the 
quasilinear system (\ref{quasi}) whose infinitesimal generators we write as
$$
  \xi(x,t) \partial_x +
  \tau(x,t) \partial_t +
  \psi(x,t,u) \partial_u
$$

Applying $\lie{}$ to
$
  \Ui_\eta  J\Xi_\eta  =
  \begin{bmatrix} \U_{\eta x} & \U_{\eta t}
  \end{bmatrix} + \U_{\eta u}
  \begin{bmatrix} u_x & u_t
  \end{bmatrix},
$
which is just the short form of (\ref{xxxx}), we obtain 
$$
  \lie{\Ui_\eta } +
  \begin{bmatrix} u_x & u_t
  \end{bmatrix} \lie J\Xi_\eta  =
  \begin{bmatrix} \psi_x & \psi_t
  \end{bmatrix} + \psi_u
  \begin{bmatrix} u_x & u_t
  \end{bmatrix}
$$ 

Thus,
\begin{eqnarray*}
  \lie{\Ui_\eta }
  \begin{bmatrix} 1\\0
  \end{bmatrix} &=& \psi_x + \psi_u u_x - \xi_x u_x - \tau_x u_t \\
  \lie{\Ui_\eta }
  \begin{bmatrix} 0\\1
  \end{bmatrix} &=& \psi_t + \psi_u u_t - \xi_t u_x - \tau_t u_t
\end{eqnarray*}

Hence by (\ref{Deltaeta}) we get the following expression for 
$\lie{\tilde \Delta_\eta }$:

\begin{longmath}
   \lie{\Ui_\eta }
  \begin{bmatrix} 0\\1
  \end{bmatrix} + \lie{(A(\U_\eta ))}
  \begin{bmatrix} u_x & u_t
  \end{bmatrix}
  \begin{bmatrix} 1\\0
  \end{bmatrix} +
  A(u) \lie{\Ui_\eta }
  \begin{bmatrix} 1\\0
  \end{bmatrix} =
  \dimension{\psi_t}s1 +
  \dimension{\psi_u}ss
  \dimension{u_t}s1 - \xi_t
  \dimension{u_x}s1 - \tau_t
  \dimension{u_t}s1 +
  \sum_i\psi^i
  \dimension{\frac
  {\partial A}
  {\partial u^i}}ss
  \dimension{u_x}s1 +
  \dimension{A(u)}ss \left(
    \dimension{\psi_x}s1 +
    \dimension{\psi_u}ss
    \dimension{u_x}s1 - \xi_x
    \dimension{u_x}s1 - \tau_x
    \dimension{u_t}s1
  \right) = \left(
    \psi_t + A(u) \psi_x
  \right) + \left(
    \psi_u - \tau_t  - \tau_x A(u)
  \right) u_t + \left(
    A(u) \psi_u - \xi_t \eigmat_s - \xi_x A(u) +
    \sum_i\psi^i \frac{\partial A}{\partial u^i}
  \right) u_x,
\end{longmath}
(where $\dimension{a}bc$ is a reminder that $a$ is $b \times c$ matrix).
Note that in this expression the coefficient of $u_t$ equals $\lie Q_\eta $.

Thus the determining system
$(\lie{\tilde \Delta})_\eta  = \lie Q_\eta \cdot\Delta$
for an infinitesimal projectable symmetry reads
\begin{eqnarray*}
  \psi_t + A(u) \psi_x &=& 0 \\
  A(u) \psi_u - \xi_x A(u) - \xi_t \eigmat_s +
  \sum_i\psi^i \frac{\partial A}{\partial u^i} &=&
  \left(
    \psi_u - \tau_t \eigmat_s - \tau_x A(u)
  \right) A(u)
\end{eqnarray*}
or
\def\bracket#1{\left[ #1 \right]}
\begin{eqnarray*}
  \psi_t + A(u) \psi_x &=& 0 \\
  \bracket{A(u),\psi_u} +
  \sum_i\psi^i \frac{\partial A}{\partial u^i} &=&
  \left(
    \xi_t \eigmat_s - \tau_t A(u)
  \right) + \left(
    \xi_x \eigmat_s - \tau_x A(u)
  \right) A(u)
\end{eqnarray*}
\et

\bex 
We continue to analyze the last determining system
in the case of a
strictly hyperbolic conservation law
\begin{eqnarray}
\label{two1} u_t + f (u,v)_x &=& 0, \\
\label{two2} v_t + g (u,v)_x &=& 0.\;\;\;
\end{eqnarray}
corresponding to (\ref{quasi}) with
$$
A = \begin{bmatrix}
   f_u & f_v \\
    g_u & g_v
  \end{bmatrix}
$$
with the characteristic values $\lambda _{1}(u,v)
< \lambda _{2}(u,v), (u,v) \in H,$ where $H$ is open in $\R^{2}$
(cf.\ \cite{Bres,Ke,KeKr}; for generalized solutions, see \cite{Ned}).

We denote by $r_{i}= {r^{1}_{i} \choose r^{2}_{i} } , \;
\ell _{i} =
{\ell ^{1}_{i} \choose \ell ^{2}_{i}}  \in ({\cal C}^{\infty } (H))^{2\times 1}$ the
characteristic vectors:
$Ar_{i} = \lambda _{i}r_{i},\;\; \ell ^{T}_{i} A =
\lambda _{i} \ell ^{T}_{i},\;  i=1,2.$

We will calculate the coefficients of an infinitesimal symmetry
$
  {\mathbf v} =
  \xi(x,t) \partial_x +
  \tau(x,t) \partial_t +
  \phi(x,t,u,v) \partial_u +
  \psi(x,t,u,v) \partial_v
$
for
(\ref{two1}),
(\ref{two2}).
By \ref{infex1} we have
\begin{eqnarray}
\label{first}
  \begin{bmatrix}
    \phi_t \\
    \psi_t
  \end{bmatrix} + A
  \begin{bmatrix}
    \phi_x \\
    \psi_x
  \end{bmatrix} &=& 0 \\
\label{second}
%  \bracket{
[A,B]
%} 
+ \phi A_u + \psi A_v &=&
  (\xi \eigmat_2 - \tau A)_t +
  (\xi \eigmat_2 - \tau A)_x A
\end{eqnarray}
where
$  B = \begin{bmatrix}
   \phi_u & \phi_v \\
    \psi_u & \psi_v
  \end{bmatrix}
$.

Differentiating (\ref{first}) with respect to $u$ and $v$ 
we have
\begin{eqnarray*}
  (B_t + \phi_x A_u + \psi_x A_v)
  \begin{bmatrix} 1 \\ 0
  \end{bmatrix} + A
  \begin{bmatrix} \phi_{xu} \\ \psi_{xu}
  \end{bmatrix} &=& 0 \\
  (B_t + \phi_x A_u + \psi_x A_v)
  \begin{bmatrix} 0 \\ 1
  \end{bmatrix} + A
  \begin{bmatrix} \phi_{xv} \\ \psi_{xv}
  \end{bmatrix} &=& 0
\end{eqnarray*}

Putting both together we obtain the matrix equation
\begin{equation}
\label{same}
  B_t + AB_x +
  \phi_x A_u +
  \psi_x A_v = 0
\end{equation}

Taking the derivative of (\ref{second}) with respect to $x$ and
subtracting from (\ref{same}) we arrive at
$$
  B_t + B_x A +
  (\xi \eigmat_2 - \tau A)_{xt} +
  (\xi \eigmat_2 - \tau A)_{xx} A = 0
$$
or in shorter form
\begin{equation}
\label{M}
  M_t + M_x A = 0 \quad
  \mbox{where} \quad
  M = B +
  (\xi _{x} \eigmat_2 - \tau  _{x}A)
\end{equation}
Thus we obtain the system (\ref{first}), (\ref{M}) which , while not equivalent to
(\ref{first}), (\ref{second}) (due to the differentiations used in deriving it), 
considerably facilitates the determination of infinitesimal symmetries.

By \cite{Tr}, Section 16, we have
$$
\begin{bmatrix} \phi(x,t,u,v) \\ \psi(x,t,u,v)
  \end{bmatrix}
 =
\sum ^{2}_{i=1} \alpha _{i}
(x-t \lambda _{i}(u,v),u,v)r_{i}(u,v),
\;\; \alpha _{i} \in {\cal C}^{\infty } (\R \times H)$$
$$M(x,t,u,v) = \sum ^{2}_{i=1} \beta _{i}
(x-t \lambda _{i} (u,v), u,v) \ell ^{T}_{i}(u,v) $$
$$\beta _{i} \in ({\cal C}^{\infty }
(\R \times H))^{2\times 1},
(x,t,u,v) \in \R^{2} \times H$$

In order to determine $\xi $ and $\tau $  as well as
$\alpha _{i} , \beta _{i}=
\begin{bmatrix} \beta^{1}_{i} \\ \beta^{2}_{i}
  \end{bmatrix}, \; i=1,2, $ we use
 $M = B+\xi _{x}I_{2} - \tau _{x}A$:
$$\sum ^{2}_{i=1} \beta ^{1}_{i} \ell ^{1} _{i} =
\sum ^{2} _{i=1} (\alpha _{i} r^{1}_{i})_{u} +
\xi _{x} - \tau _{x} f_{u}; \;
\sum ^{2}_{i=1} \beta ^{1}_{i} \ell ^{2}_{i} =
\sum ^{2} _{i=1} (\alpha _{i} r^{1}_{i}) _{v}
-\tau _{x} f_{v}$$
$$\sum ^{2}_{i=1} \beta ^{2}_{i} \ell ^{1}_{i} =
\sum ^{2} _{i=1} (\alpha _{i} r^{2}_{i}) _{u}
-\tau _{x} g_{u};\;
\sum ^{2}_{i=1} \beta ^{2}_{i} \ell ^{2}_{i} =
\sum ^{2} _{i=1} (\alpha _{i} r^{2}_{i}) _{v}
+\xi _{x} -\tau _{x} g_{v}  $$

By \ref{factorinfprop} resp.\ \ref{factorinfpropG}, solutions
of this system in ${\cal C}^\infty$ resp.\ $\G$ determine infinitesimal
symmetries in each of the respective settings.
\et

Note that  the above  considerations also hold for
systems of order $s>2$ even in non-conservative
form with the additional assumption that the coefficients
of $A$ do not depend on $x$ and $t$.

\end{document}